\documentclass[12pt,twoside]{amsart}
\usepackage{amssymb}

\numberwithin{equation}{section}
\font\tengothic=eufm10 scaled\magstep 1
\font\sevengothic=eufm7 scaled\magstep 1
\newfam\gothicfam
\textfont\gothicfam=\tengothic
\scriptfont\gothicfam=\sevengothic

\newtheorem{theorem}{Theorem}[section]
\newtheorem{lemma}[theorem]{Lemma}
\newtheorem{proposition}[theorem]{Proposition}
\newtheorem{corollary}[theorem]{Corollary}
\newtheorem{conjecture}[theorem]{Conjecture}

\theoremstyle{definition}
\newtheorem{definition}[theorem]{Definition} 
\newtheorem{remark}[theorem]{Remark}
\newtheorem{example}[theorem]{Example}

\newtheorem{notation}[theorem]{Notation}

\newcommand{\reg}{\operatorname{reg}}

\newcommand{\hu}{\underline h}

\newcommand{\Tor}{\operatorname{Tor}}

\newcommand {\PP}{\mathbb{P}}

\def\T{{\it T}}
\def\P{{\mathbb P}}

\def\V{{\mathbb V}}
\def\X{{\mathbb X}}

\def\Y{{\mathbb Y}}

\def\Z{{\mathbb Z}}

\def\u{\underline}

\begin{document}
\title[First infinitesimal neighborhood]{On the first infinitesimal 
neighborhood of a linear configuration of points in $\P^2$ }

\author[A.V.\ Geramita, J.\ Migliore, L.\ Sabourin]{A.V.\ Geramita, J.\
Migliore, L.\ Sabourin}

\thanks{The first author was supported, in part, by grants from NSERC 
(Canada) and INDAM (Italy). ÊPart of the work for this Êpaper was done 
while the second author was sponsored by the National Security Agency (USA) 
under Grant Number MDA904-03-1-0071. The third author was supported, during 
the writing of this paper, by an NSERC (Canada) PDF at Notre Dame and York 
Universities.}


\begin{abstract}
We consider the following open questions. ÊFix a Hilbert function $\hu$, 
that occurs for a reduced zero-dimensional subscheme of $\P^2$. ÊAmong all 
subschemes, $\X$, with Hilbert function $\hu$, what are the possible 
Hilbert functions and graded Betti numbers for the first infinitesimal 
neighborhood, $\Z$, of $\X$ (i.e. the double point scheme supported on 
$\X$)? ÊIs there a minimum ($\hu^{\min}$) and maximum ($\hu^{\max})$ such 
function?

The numerical information encoded in $\hu$ translates to a {\it type 
vector} which allows us to find unions of points on lines, called {\it 
linear configurations}, with Hilbert function $\hu$. ÊWe give necessary and 
sufficient conditions for the Hilbert function and graded Betti numbers of 
the first infinitesimal neighborhoods of {\bf all} such linear 
configurations to be the same. ÊEven for those $\hu$ for which the Hilbert 
functions or graded Betti numbers of the resulting double point schemes are 
not uniquely determined, we give one (depending only on $\hu$) that does 
occur. ÊWe prove the existence of $\hu^{\max}$, in general, and discuss 
$\hu^{\min}$. ÊOur methods include liaison techniques.
\end{abstract}


\maketitle

\tableofcontents


ÊÊ \section{Introduction} \label{intro}

The classification of all the possible Hilbert functions of reduced
zero-dimensional subschemes of the projective space $\PP ^n(k)$, ($k$ a
field of characteristic zero) is well known (see e.g.\ \cite{GMR}). ÊIn marked
contrast, the analogous classification, even for the important class of 
non-reduced zero-dimensional schemes which are unions of ``2-fat point" 
schemes (which we will refer to as {\it double point schemes}; see $\S 2$ 
for the definitions) is wide open. This in spite of the fact that answers 
to such questions have interesting implications in algebraic geometry, 
coding theory, computational complexity and statistics.

The one major achievement in this area is the proof of the genericity
(apart from some well understood and well known exceptions) of the
Hilbert function of a generic double point scheme by J. Alexander and A.
Hirschowitz (see \cite{AH3} and also K. Chandler \cite{ chandler}). ÊOther 
important contributions to our understanding of the Hilbert function of fat 
point
schemes have been made by Bocci \cite{bocci}, Catalisano, Geramita and 
Gimigliano \cite{cgg}, Ciliberto \cite{ciliberto}, Ciliberto and Miranda 
\cite{ciliberto-miranda}, \cite{ciliberto-miranda2}, Laface and Ugaglia 
\cite{laface-ugaglia}, \cite{laface-ugaglia2}, Yang \cite{yang} (for 
further references see the survey article of Miranda \cite{miranda}).

Since one knows the functions which can be the Hilbert function of a
reduced zero-dimensional subscheme of $\PP ^n$, one may then inquire as
to the possible minimal free resolutions for reduced subschemes sharing
the same Hilbert function. ÊThis is, in general, a very difficult problem
which has attracted a great deal of attention. ÊThere are: complete
results for reduced subschemes of $\PP ^2$ by Campanella \cite{campanella}; 
a sharp upper bound for any Hilbert function (in terms of the graded Betti 
numbers) by Bigatti \cite{bigatti}, Hulett \cite{hulett} and Pardue 
\cite{pardue}; complete results in low codimension and under the assumption 
that the coordinate ring of the reduced scheme is (in some way) special -- 
e.g. for codimension two and codimension 3 Gorenstein see e.g. Diesel 
\cite{diesel} and Geramita and Migliore \cite{GM}, while for Gorenstein 
rings with the Weak Lefschetz Property see, e.g., ÊGeramita, Harima, and 
Shin \cite{GHSAdv}, Migliore and Nagel \cite{MN2}.

The results about minimal free resolutions for fat point schemes are much
scantier. ÊNotable results (in $\PP ^2$) are given by: Catalisano
(arbitrary fat points schemes supported on a conic) \cite {catalisano}; 
Harbourne et.al. (for fat point schemes supported on generic sets of points 
or arbitrary fat point schemes supported on few points on a cubic) 
\cite{harbourne}, \cite{FHH}, \cite{HHF}. ÊFor higher dimensional spaces we 
have: ÊCatalisano, Ellia and Gimigliano \cite{ cat-el-gim} \ Ê(arbitrary 
fat points whose support is on a rational normal curve in $\P^n$); Fatabbi 
\cite{fat}, Francisco \cite{fran}, Fatabbi and Lorenzini \cite{fat-lor} and 
Valla \cite {valla}\ (arbitrary fat points whose support is on $\leq n$ 
points in $\P^n$).

In this paper we will deal with the problem of classification of Hilbert
functions and resolutions of double point schemes in $\PP ^2$ in the
following way. ÊLet $\hu$ be the Hilbert function of a reduced
zero-dimensional subscheme of $\PP ^2$. ÊWhat are all the possible
Hilbert functions and minimal free resolutions for double point schemes
whose support, $\X$, has Hilbert function $\hu$?

Given any such Hilbert function $\underline{h}$, there is a well known 
family of reduced subschemes of $\P^2$ whose Hilbert function is 
$\underline{h}$ -- the so-called $k$-configurations of a specific type (see 
$\S 2$ for the definitions). ÊSo, the first natural problem to consider is 
the following:

\begin{quote}
\it{ if $\X$ is a $k$-configuration in $\P^2$ with Hilbert function $\hu$, 
can we describe the Hilbert function of the double point scheme whose 
support is $\X$?}
\end{quote}

As is also well known, the $k$-configurations in $\P^2$ of the same type 
always have the same graded Betti numbers in their minimal free resolution 
(see \cite{GHSAdv} for this and generalizations to $\P^n$). Ê This nice 
result for reduced subschemes of $\P^2$ leads naturally to another question:

\begin{quote}
{\it if $\X$ is a $k$-configuration in $\P^2$ with Hilbert function $\hu$, 
can we describe the graded Betti numbers in the minimal free resolution of 
the double point scheme whose support is $\X$?}
\end{quote}

It is also well known that all reduced sets of points which form
$k$-configurations with Hilbert function $\hu$ enjoy many extremal
properties (e.g. their minimal free resolutions are the extremal ones
described by Bigatti, Hulett and Pardue (see \cite{GHSAdv}), their
conductors are extremal - at least in $\PP ^2$ (see \cite{bazz})). ÊIt is
thus an obvious question to ask if the double point schemes supported on
linear configurations have any extremal properties.

Indeed, having raised the problem of extremal properties, one is naturally 
led to another series of questions:

\begin{quote}
{\it if $\hu$ is the Hilbert function of some reduced zero-dimensional 
subscheme of $\P^2$, do there exist Hilbert functions $\u{h}^{\max}$ and 
$\u{h}^{\min}$ of double point schemes (whose support has Hilbert function 
$\u{h}$) such that if ${h}^\prime$ is the Hilbert function of any double 
point scheme whose support has Hilbert function $\u{h}$, then
$$
\u{h}^{\min} \leq h^\prime \leq \u{h}^{\max}.
$$
Moreover, if $\u{h}^{\min}$ and $\u{h}^{\max}$ exist, what are they? Ê}
\end{quote}

The questions above are the main ones we will consider in this paper.

\medskip

We now give a summary description of the results we have obtained 
below. ÊThe first thing we do is restrict our notion of a $k$-configuration 
to that of a {\it linear configuration} (see Definition 
\ref{k-config}). ÊThis technical restriction is crucial for our main 
theorems and yet does not restrict the Hilbert functions we consider for 
reduced sets of points.

If $\u{h}$ is the Hilbert function of a zero dimensional reduced subscheme 
of $\P^2$ we first describe an O-sequence (see Definition \ref{standard 
os}), $dbl(\u{h})$, with the property that: if $\X$ is a special linear 
configuration with Hilbert function $\u{h}$ (and $\X$ always exists) then 
the double point scheme with support $\X$ has Hilbert function $dbl(\u{h})$ 
(see Theorem \ref{any 2-type}).

We then give a complete description of all the Hilbert functions $\u{h}$ so 
that {\it every} double point scheme whose support is a linear 
configuration with Hilbert function $\u{h}$, has Hilbert function 
$dbl(\u{h})$ (see Theorem \ref{fat = bdl}).

Even when it is no longer true that every double point scheme whose support 
is a linear configuration with Hilbert function $\hu$ has the same Hilbert 
function, we prove that all such double point schemes share the same 
regularity (see Theorem \ref{same reg}), which is the maximal possible for 
double point schemes whose support has Hilbert function $\hu$. ÊThis 
illustrates one sort of ``extremality" property for the Hilbert functions 
of double point schemes whose support share the same Hilbert function.

We also investigate the minimal free resolutions of double point schemes 
supported on a linear configuration. ÊWe give necessary and sufficient 
conditions on the Hilbert function $\hu$ in order that all double point 
schemes with support on a linear configuration having Hilbert function 
$\hu$ have the same graded Betti numbers in their minimal free resolution 
(see Theorem \ref{fat = bdl}). ÊAs expected, the results about minimal free 
resolutions are more subtle and restrictive than those simply about Hilbert 
functions.

It is worth making some comment here about our method of proof for the 
results about double point schemes sharing the same Hilbert function for 
their support and, in particular, when the support is a linear configuration.

Although our principal aim in this paper is the study of the possible 
Hilbert functions of double point schemes in $\P^2$, we spend a great deal 
of effort (especially in $\S 3$ and $\S 4$) studying special configurations 
of {\bf reduced} point sets in $\P^2$ (which we call {\it pseudo linear 
configurations}). ÊAlthough these reduced point sets are really peripheral 
to our main concern, there are important reasons for considering them which 
come out of the strong connections between the numerical information 
encoded in these reduced schemes and the numerical information we seek 
about the 2-fat point schemes we
are considering. ÊIn fact, our approach illustrates (in a very concrete 
way) how one can get a great deal of mileage out of considering certain 
collections of Ê2-fat point schemes in $\P^2$ as if they were the union of 
a collection of triples of reduced points (configured in a special 
way). ÊThis sort of philosophy is evident in J. Alexander and A. 
Hirschowitz's ``{\it Horace Method Ê(Divide and Conquer)} " \cite{AH3} and 
also in \cite{cgg2} and \cite{dent}. ÊThe novelty of our approach is that, 
for the first time, we use the techniques of Liaison as an additional 
weapon for {\it Horace's} arsenal.

In $\S 7$ we consider the problem of existence for, what we have called, 
${\hu}^{\max}$ and ${\hu}^{\min}$. ÊWe note that ${\hu}^{\max}$ always 
exists, even though it is difficult in practice to say exactly what it 
is. ÊRecall that the results of J. Alexander and A. Hirschowitz (see 
\cite{AH3}) give (as a special case in $\P^2$) that: if we denote by 
$\hu_s$ the generic Hilbert function of a set of $s$ distinct points in 
$\P^2$, then $\hu_s^{\max} = \hu_{3s}$ Êexcept for $s = 2,\ 5$. ÊWe have 
been unable to decide if ${\hu}^{\min}$ exists, in general. ÊNevertheless, 
we have found it for $\hu_s$ when $s = \binom{t}{2}$. ÊEven though we 
cannot decide if ${\hu}^{\min}$ actually exists, in general, we can prove 
something that would be a consequence of that existence: namely the 
existence of a maximal regularity for all double point schemes whose 
support has Hilbert function $\hu$ (see Remark \ref{same reg}).

\section{Preliminaries} \label{prelim section}

Let $k$ be any infinite field of characteristic zero and let $R =
k[x_0,x_1,x_2]$. ÊWe denote by $\PP ^2(k)$ the scheme $proj(R)$. ÊIf $P$ is a
point in $\PP^2$ defined by the prime ideal $\wp = (L_1, L_2)$ (the $L_i$
linearly independent linear forms in $R$) then any scheme supported on the
point $P$ is defined by a $\wp$-primary ideal of $R$.

\begin{definition}
A scheme supported on the point $P$ is called a {\em fat point scheme with
support $P$} if it is defined by the primary ideal
$\wp ^t$ for some integer
$t > 1$. If, in particular, $t = 2$ then we shall call the scheme defined by
$\wp^2$ a {\em double point scheme with support on $P$}. ÊThis latter is also
referred to as the {\em first infinitesimal neighborhood of $P$}.

More generally, if $\X = \{ P_1, \ldots , P_\ell \}$ is any set of distinct
points in $\PP ^2$ where $P_i$ is defined by the prime ideal $\wp_i$, then the
{\em double point scheme with support $\X$} is the scheme defined by the
(saturated) ideal $\wp_1^2 \cap \ldots \cap \wp_\ell^2$.

If a scheme supported on $\X$ is defined by an ideal of the type $\wp_1^{n_1}
\cap \ldots \cap \wp_\ell^{n_\ell}$ then we sometimes loosely refer to it as a
{\em fat point scheme} with support $\X$. ÊIf, in addition, the $n_i $ are 
all
the same (and say are equal to $t$) then we say that the scheme defined on $\X$
is a {\em $t$-fat point scheme on $\X$}.
\end {definition}

We also recall some terminology that is used in discussing the Hilbert 
function
of zero dimensional subschemes of $\PP ^2$.

\begin{definition}\label{hilb} Let $\hu$ be the Hilbert function of a zero 
dimensional subscheme, say $\X$, of $\PP^2$. ÊWe define:
\begin{itemize}
\item[({\it i})] $\alpha (\hu)$ to be the least integer\ $t$ \ for which
$\hu(t) < \binom{t+2}{2}$;

\item[({\it ii})] $\Delta\hu$ to be the {\it first difference} of $\hu$, i.e.
$$
\Delta\hu(t) = \hu(t) - \hu(t-1);
$$

\item[({\it iii})] $\sigma(\hu)$ to be the least integer $t$
for which $\Delta\hu(t) = 0$.

\end{itemize}

We also sometimes refer to $\hu$ as $h_{\X}$. ÊIn this case, since $\Delta 
h_{\X}(t)\neq 0$ for only finitely many values of $t$, we refer to the sequence
$$
\Delta h_{\X}(0) = 1 \ \ \ \ \Delta h_{\X}(1) \ \ \ \ \cdots \ \ \ \ \Delta 
h_{\X}(\sigma(h_{\X}) -1)
$$
as the {\it $h$-vector of $\X$}.

If the scheme $\X$ is defined by the ideal $I$ of the ring $R$ we will also 
use the notation
$h_{\X} = h_{R/I}$.

\end{definition}

Geramita, Harima and Shin defined the notion of an $n$-type vector in
\cite{ghs}. ÊSince we only need the case of a 2-type vector, we only recall
that definition.

\begin{definition}
A {\em 2-type vector} is a vector of the form $T=(d_1,d_2,\ldots ,d_r)$, where
$0<d_1<d_2<\ldots <d_r$ are integers. ÊFor such a 2-type vector, we define 
$\alpha (T)=r$ and $\sigma (T)=d_r$.
\end{definition}

\begin{theorem}\label{1-1 co} \cite[Theorem 2.6]{ghs} Let $S_2$ denote the
collection of Hilbert functions of all sets of distinct points in $\PP ^2$.
Then there is a 1-1 correspondence $S_2\leftrightarrow \{ 2\mbox{-type
vectors}\}$. ÊUnder this correspondence if $\hu\in S_2$ and $\hu$ 
corresponds to
$T$ (we write $\hu\leftrightarrow T$) then $\alpha (\hu) =\alpha (T )$ and 
$\sigma
(\hu) = \sigma (T )$.
\end{theorem}

We now give the inductive formula for obtaining the Hilbert function referred
to Êin Theorem \ref{1-1 co} from its corresponding $2$-type vector.

\begin{theorem}\cite[Proof of Theorem 2.6]{ghs}\label{formula}
Let $T =(d_1,d_2,\ldots d_r)$ be a 2-type vector, and let $\hu_i$ denote the
Hilbert function of $d_i$ points on a line. ÊThen $\hu\leftrightarrow T$ where
$\hu(j)=\hu_r(j)+\hu_{r-1}(j-1)+\ldots +\hu_1(j-(r-1))$ and $\hu(t)=0$ for 
$t<0$.
\end{theorem}

The notion of an $n$-type vector is convenient for defining the notion of a
$k$-configuration in $\PP^n$. ÊWe give here the definition of a
$k$-configuration in $\PP^2$.

\begin{definition}\label{k-config}

\noindent $a)$\ Let $T=(d_1,d_2,\ldots ,d_r)$ be a 2-type vector. ÊLet 
$L_1,L_2, \ldots ,L_r$ be Êdistinct lines in $\PP^2$. ÊLet $\X_i$ consist 
of $d_i$ points on $L_i$ for each $i$. ÊSuppose, furthermore, that $L_i$ 
does not contain any point of $\X_j$ for $j<i$. ÊThen $\X=\cup_{i=1}^r 
\X_i$ is called a {\em $k$-configuration of type $T$}.

\noindent $b)$\ If we assume further that no point of $\X_i$ is on line 
$L_j$, for $j \neq i$, then $\X$ will be called a {\it linear configuration 
of type $T$}.
\end{definition}

\begin{example}\label{order}
In the diagram below, $\X _1$ consists of the two points of $L_1$ that are not
on $L_2$, $\X _2$ consists of the five points of $L_2$, and $\X _3$ consists of
the 6 points of $L_3$. ÊThen $\X =\X _1 \cup \X _2 \cup \X _3$ is a
$k$-configuration of type $\T =(2,5,6)$. ÊNotice that $L _i$ does not contain a
point of $\X _j $ for $j<i$, although $L _1$ does contain a point of $\X 
_2$. ÊThus $\X$ is {\bf NOT} a linear configuration of type $T = (2,5,6)$.

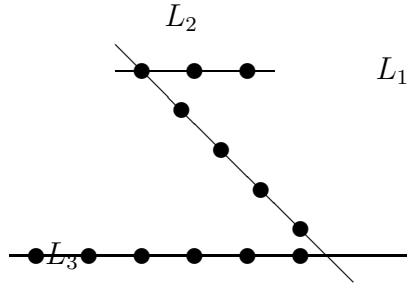
\vbox{$$
\begin{picture}(130,140)(-20,-100)
\put(-10,0){\line(1,0){60}}
\put(-10,10){\line(1,-1){90}}
\put(-50,-70){\line(1,0){150}}
\put(0,0){\circle*{6}}
\put(20,0){\circle*{6}}
\put(40,0){\circle*{6}}
\put(15,-15){\circle*{6}}
\put(30,-30){\circle*{6}}
\put(45,-45){\circle*{6}}
\put(60,-60){\circle*{6}}
\put(-40,-70){\circle*{6}}
\put(-20,-70){\circle*{6}}
\put(0,-70){\circle*{6}}
\put(20,-70){\circle*{6}}
\put(40,-70){\circle*{6}}
\put(60,-70){\circle*{6}}
\put(65,-5){\makebox(0,0)[bc]
ÊÊ Ê Ê Ê Ê Ê Ê{\makebox(16,12)[bl]{{$L_1$}}}}
\put(-15,15){\makebox(0,0)[bc]
ÊÊ Ê Ê Ê Ê Ê Ê{\makebox(16,12)[bl]{{$L_2$}}}}
\put(-60,-75){\makebox(0,0)[bc]
ÊÊ Ê Ê Ê Ê Ê Ê{\makebox(16,12)[bl]{{$L_3$}}}}
\end{picture}$$}

\noindent
Notice that $\X$ is not a $k$-configuration of type $\T$ = (3,4,6) since Ê$\X
_1$ would have to consist of all 3 points on $L _1$ and this includes a point
of $L _2$. Ê \qed
\end{example}

One can see from Example ~\ref{order} that the fact that the $\X _i$ are
ordered from smallest to largest is a crucial part of the definition of a
$k$-configuration. ÊAs well, the example suggests that the same
$k$-configuration cannot have two different 2-type vectors associated to it.
In fact, more is true: namely, all $k$-configurations of type $T$ have the
Hilbert function corresponding to $T$.

\begin{theorem}[\cite{ghs}, p. 21] \label{hilbkconfig}
Let $\X$ be a $k$-configuration of type $T\leftrightarrow \hu$. ÊThen 
$h_{\X}=\hu$.
\end{theorem}

\begin{remark} Although some of the results of this paper (and results 
cited from earlier papers) are true for arbitrary $k$-configurations, our 
main results are not. ÊFor this reason, {\bf from now on we will only 
consider linear configurations} (see Definition \ref{k-config}).
\end{remark}

Recall that the $(i,j)^{th}$ graded Betti number of an ideal $I$ of $R$ is
defined to be
$$
\beta_{i,j}^I :=(\Tor_i(R/I,k))_j.
$$
entry 

It turns out that the graded Betti numbers of a linear configuration of 
type $T$
are also completely determined by $T$ (\cite{GHSAdv}, Theorem 3.6). ÊIn fact,
those Betti numbers are extremal, in a way which we will explain later. ÊWe
will see in this paper exactly when the Hilbert function (Corollary ~\ref{fat =
bdl} (a)) and graded Betti numbers (Corollary ~\ref{fat = bdl} (b),(c)) of the
double points supported on a linear configuration are determined just from the
type of the linear configuration - something that does not always happen!

Even when the Hilbert function of double points supported on a linear
configuration is not determined simply from the type of the linear
configuration, we will at least be able to determine the Hilbert function
of double points supported on very special linear configurations. ÊWe now 
proceed
to the definitions of these two special classes of linear configurations (see
\cite[before Example 4.1]{GS}).

\begin{definition}\label{std k-config def}
ÊÊA linear configuration of type $T=(d_1,d_2,\ldots ,d_r)$ in $\PP^2$ is 
called a {\em standard linear configuration of type $T$} if it consists of:
$$
\begin{array}{llcl}
d_r& \mbox{ points with coordinates }&[j:0:1]& 0\leq j\leq d_r -1, j \in
\mathbb N,\\
\\
ÊÊ Ê Ê & Ê Ê Ê& \vdots Ê Ê Ê Ê Ê Ê Ê Ê Ê Ê Ê Ê Ê Ê Ê Ê Ê Ê& Ê\\ \\
d_2& \mbox{ points with coordinates }&[j:r-2:1]& 0\leq j\leq d_2 -1, j \in
\mathbb N,\\ \\
d_1& \mbox{ points with coordinates }&[j:r-1:1]& 0\leq j\leq
d_1 -1, j \in \mathbb N.
\end{array}
$$
\end{definition}

\begin{definition} Let $J$ be a homogeneous ideal in $S=k[x_1,\ldots ,x_n]$.
We say that a radical ideal $I$ of $R=k[x_0,x_1,\ldots ,x_n]$ {\em lifts} $J$
if there is a linear form $L$ which is a non-zero-divisor on $R/I$ for which
$(I,L)/L\simeq J$.
\end{definition}

If $I$ is an ideal of $R=k[x_0,\ldots ,x_n]$ which lifts the homogeneous ideal
$J$ of $S=k[x_1,\ldots ,x_n]$, then the minimal free $R$-resolution of $I$ has
the same graded Betti numbers as the minimal free $S$-resolution of $J$ (see
\cite{bruns-herzog}, Proposition 1.1.5).

Note that the ideal of the standard linear configuration of type $T=(d_1,d_2,
\ldots ,d_r)$ is a lifting of the monomial ideal $\langle
x^{d_r},x^{d_{r-1}}y,x^{d_{r-2}}y^2,\ldots ,y^{r}\rangle$ (an ideal is called
{\em monomial} if it is generated by monomials). ÊWe call this the {\em
standard lifting}.

Note also that the monomial ideal being lifted to obtain the standard
linear configuration is by no means random, but rather satisfies the following
very special condition: if a monomial $m\in I$, then every larger monomial
(using the lexicographic ordering) of the same degree is also in $I$. ÊSuch
ideals are called {\em lex-segment ideals}.

Since the ideal of a standard linear configuration always lifts a lex-segment
ideal (by \cite{GHSAdv}, Theorem 4.3), standard linear configurations can be
looked at as providing the 1-1 correspondence between Hilbert functions of
points and lex-segment ideals.

The special linear configurations for which we will always be able to 
determine
the Hilbert functions of the double points with that support are defined in
almost the same way as standard linear configurations.

\begin{definition} \label{spread out config} ÊA linear configuration of type
$T=(d_1,d_2,\ldots ,d_r)$ in $\PP^2$ is called a {\em spread out linear 
configuration of type $T$} if it consists of:
$$
\begin{array}{llcl}
d_r& \mbox{ points with coordinates }&[j:0=d_r-d_r:1]& 0\leq j\leq d_r -1, j
\in \mathbb N,\\ \\
ÊÊ Ê Ê & Ê Ê Ê & \vdots Ê Ê Ê Ê Ê Ê Ê Ê Ê Ê Ê Ê Ê Ê Ê Ê Ê Ê& Ê\\
d_2& \mbox{ points with coordinates }&[j:d_r-d_2:1]& 0\leq j\leq d_2 -1, j \in
\mathbb N,\\ \\
d_1& \mbox{ points with coordinates }&[j:d_r-d_1:1]& 0\leq
j\leq d_1 -1, j \in \mathbb N.
\end{array}
$$
\end{definition}

\begin{example} If $T =(1,2,4,7)$, then the standard linear configuration 
and the spread out linear configuration of type $T$ are as follows:
$$\begin{array}{cc}
\mbox{\underline{standard}}&\mbox{\underline{spread out}}\\
ÊÊ& \\
\begin{array}{ccccccc}
ÊÊ& & & & & & \\
ÊÊ& & & & & & \\
ÊÊ& & & & & & \\
\bullet& & & & & & \\
\bullet&\bullet& & & & & \\
\bullet&\bullet&\bullet&\bullet& & & \\
\bullet&\bullet&\bullet&\bullet&\bullet&\bullet&\bullet
\end{array} \hspace{.5cm} Ê&\hspace{.5cm}
\begin{array}{ccccccc}
\bullet& & & & & & \\
\bullet&\bullet& & & & & \\
\circ&\circ&\circ& & & & \\
\bullet&\bullet&\bullet&\bullet& & & \\
\circ&\circ&\circ&\circ&\circ& & \\
\circ&\circ&\circ&\circ&\circ&\circ& \\
\bullet&\bullet&\bullet&\bullet&\bullet&\bullet&\bullet
\end{array}
\end{array}
$$
where the $\circ$'s represent ``imaginary" points that we are using to
properly Êposition the points in which we are interested. ÊAgain in this case,
we want rows consisting of 1, 2, 4 and 7 points but to obtain the spread out
linear configuration of this type we add rows of 3, 5 and 6 ``imaginary" 
points to
form an ``isosceles right triangle". \qed
\end{example}

\begin{remark}
The process of forming a spread out linear configuration ensures that the
``diagonal" Êpoints are collinear and it is this fact that will enable us to
determine the Hilbert function of sets of double points with support on a
spread out linear configuration. \qed
\end{remark}

The notion of {\em basic double linkage} is extremely useful, both in liaison
theory (where it is fundamental) and as a construction tool for building
interesting schemes. ÊWe use it in this paper to construct sets of double
points. ÊBecause we are primarily interested in points in $\mathbb P^2$, we
recall the basic ideas here only in that context, even though they are
applicable in far greater generality (cf.\ \cite{LR}, \cite{mig-book},
\cite{KMMNP}, \cite{BM4}, \cite{GM4} for basic properties). ÊWe collect the
known facts about basic double linkages in ${\mathbb P}^2$ here for the
convenience of the reader. ÊIf no other reference is given, see \cite{mig-book}
for details.

\begin{lemma}[Basic Double Linkage]\label{BDL lemma}

Let $\X$ be a zero dimensional subscheme of $\mathbb P^2$. ÊLet $F \in I_\X$
be any polynomial, and let $G$ be any polynomial such that $\{ F,G \}$ form a
regular sequence. Ê(It makes no difference if $G$ vanishes on a point of $\X$
or not.) ÊForm the ideal $I = G \cdot I_\X + (F)$. ÊThen

\begin{itemize}
\item[(a)] $I = I_\Z$ is the saturated ideal of a subscheme $\Z$ in
$\mathbb P^2$.

\item[(b)] The support of $\Z$ is the union of the support of $\X$ and the
support of the complete intersection scheme, $\V$, defined by $(F,G)$.

\item[(c)] If $\deg F = d_1$ and $\deg G = d_2$ then there is an exact
sequence
\begin{equation} \label{bdl exact seq}
0 \rightarrow R(-d_1-d_2) \rightarrow I_\X(-d_2) \oplus R(-d_1) \rightarrow
I_\Z
\rightarrow 0
\end{equation}

\item[(d)] We have the Hilbert function formula
\begin{equation} \label{hf formula}
h_{R/I_\Z}(t) = h_{R/(F,G)}(t) + h_{R/I_\X}(t-d_2).
\end{equation}
(We will often use the first difference of this formula which, for example,
gives that $\deg \Z = \deg \X + d_1d_2$.)

\item[(e)] Ê(\cite{MN1} Corollary 4.5) Suppose that $I_\X$ has a minimal free
resolution
$$
0 \rightarrow{\mathbb F}_2\rightarrow{\mathbb F}_1\rightarrow 
I_\X\rightarrow 0.
$$
Then $I_\Z$ has a free resolution
$$
0\rightarrow R(-d_1-d_2)\oplus{\mathbb F}_2(-d_2)\rightarrow{\mathbb F}_1(-d_2)
\oplus R(-d_1)\rightarrow I_\Z\rightarrow 0.
$$
Furthermore, this resolution is minimal if and only if $F$ is {\em not} a
minimal generator of $I_\X$. ÊIf $F$ is a minimal generator of $I_\X$ then one
term, $R(-d_1-d_2)$, splits off, yielding a minimal free resolution.
\end{itemize}
\end{lemma}

\begin{remark}
ÊÊIt is easy to see that any standard linear configuration or any spread out
linear configuration, $\X$, can be produced by a sequence of basic double
linkages. ÊSimply start at the ``top" and choose as the polynomials $F$
suitable unions of ``vertical" lines, and choose as $G$ sequentially the
``horizontal" lines containing points of $\X$, working down from top to
bottom. Ê(See Proposition \ref{std pseudo has std os} below.)

Less obviously, as was observed by one of us (Migliore) several years ago, 
{\em any} linear configuration in $\mathbb P^2$ can be viewed as a sequence 
of basic double links. That fact Êwill be seen later as a consequence of 
the more general result in Theorem \ref{pseudo k con = bdl}.

But, the main new idea in this paper comes out of the realization that many
double point schemes can also be obtained as the result of a sequence of basic
double links. ÊSince this idea is pervasive in this paper, it will be useful to
have a simple example that illustrates the point. Ê\qed
\end{remark}

\begin{example}\label{build fat}
We construct a 2-fat point scheme whose underlying support is a linear
configuration of type $(1,2)$. ÊWe shall do this example in some detail as
it illustrates, in a simple way, some of the key ideas of the proofs in
this paper. ÊIn particular, it illustrates how: basic double links can be used
to ``fatten up" points and how one can use basic double links (with the form
$F$ progressively growing in degree) to get linear configurations of 2-fat 
points.

Without loss of generality we may assume that our points are
\[
P_1 = [1:0:0], \ \ \ P_2 = [0:1:0], \ \ \ P_3 = [0:0:1].
\]
We want to ``fatten up" $P_1$ by ``adding'' to it something of length 2 and 
we then
want to fatten up $P_2$ and $P_3$ by adjoining to each a length 2 piece.

So, it is as if we were considering $1 + 2$ points on the first 
``horizontal" line and
then $2+4$ points on the second line. ÊWe write those numbers down, 
$1,2,2,4$, and use
them as a guide for our construction.

Using the ``1'' we begin with the point $P_1 = [1:0:0]$, and the ideal 
$I_{P_1} =
(y,z)$. ÊNow consider the 2's. ÊWe take them both and think of performing a 
basic
double link on $I_{P_1}$ which will, at the same time, ``fatten up" $P_1$ 
(the first 2
in our sequence) and add the two reduced points $P_2$ and $P_3$ on the line 
$x$ (the
second 2 in the sequence).

Let $F = yz$ and $G = (y+z)x$ and form $I = G\cdot I_{P_1} + (F)$.
As noted in Lemma \ref{BDL lemma}, $I$ is the saturated ideal of a scheme 
supported on the union of $P_1$ and the support of the Êscheme defined by 
$(F,G)$. ÊThe latter scheme is supported on
$P_1,P_2,P_3$.

Now,
\[
\begin{array}{rcl}
I & = & ((y+z)xy, (y+z)xz, yz) \\
& = & (y^2x, z^2x, yz) \\
& = & I_{P_1}^2 \cap I_{P_2} \cap I_{P_3}.
\end{array}
\]
The last equality can be checked locally, since we know that $I$ is 
saturated. ÊNow we
use the ``4'' to fatten up $P_2$ and $P_3$, by letting $F = yz(x+z)(x+y), G 
= x$.
Clearly $F \in I$ and we use $F$ and $G$ to form a basic double link on 
$I$. ÊWe obtain
\[
\begin{array}{rcl}
J & = & xI +(F) \\
& = & (y^2x^2, z^2x^2, xyz, yz(x+z)(x+y)) \\
& = & (y^2x^2, z^2x^2, xyz, y^2z^2) \\
& = & I_{P_1}^2 \cap I_{P_2}^2 \cap I_{P_3}^2
\end{array}
\]
as we wanted to show. \qed

\end{example}

\begin{remark} Finally, we recall that if $\Z \subset {\mathbb P}^2$ then 
$I_\Z$ has
{\em regularity}
$d$ if
\[
d = \min \{ t \ | \ h^1({\mathcal I}_\Z (t-1)) = 0 \}.
\]
If this is the case then $I_\Z$ is generated in degrees $\leq d$ 
(\cite{mumford}).
Furthermore,
\[
d = \min \{ t \ | \ \Delta h_{R/I_\Z}(t) = 0 \} = \sigma(h_{R/I_{\Z}}).
\]
We will say that $\Z$ has regularity $d$ if $I_\Z$ does. Ê\qed
\end{remark}

ÊÊThe following elementary result about the regularity of the first 
infinitesimal neighborhood of a set of distinct points in $\P^2$ will be 
extremely useful.

\begin{lemma} \label{bd reg}
Let $\X$ be a reduced set of points in $\mathbb P^2$ with regularity
$r+1$. ÊLet $\Z$ be the first infinitesimal neighborhood of $\X$. ÊThen
$\reg(I_{\Z}) \leq 2 \cdot \reg(I_\X) = 2r+2$.
\end{lemma}

\begin{proof}
By hypothesis, $\X$ imposes independent conditions on forms of degree
$r$. ÊWe want to show that $\Z$ imposes independent conditions on forms
of degree $2r+1$. ÊThis means that we want to show that if $P \in \X$ and
$\Z'$ is the subscheme of $\Z$ supported on $\X' = \X \backslash P$ then
there is a form of degree $2r+1$ vanishing on $\Z'$ and also on $P$
together with any tangent direction at $P$. ÊBut this is clear: let $F$
be a form of degree $r$ vanishing on $\X'$ but not at $P$. ÊThe $F^2$
vanishes on $\Z'$ but not at $P$, and if $L$ is the line through $P$ with
the desired tangent direction then $F^2L$ is the desired form.
\end{proof}

\vskip 1cm


\section{Pseudo linear configurations}

Before we can begin to consider configurations of 2-fat points in the 
plane, it
is useful for us to extend the class of configurations of simple points in the
plane whose Hilbert function we can control. ÊThese will play an important part
in our attempt to discover the Hilbert function of all 2-fat point schemes
whose underlying supports have the same Hilbert function. ÊThe new 
configurations we consider are inspired by Example \ref{build fat}.

\begin{definition} \label{def of pseudo}
$ $

\begin{itemize}
\item[ $i)$ ] A {\em pseudo type vector} is a sequence $T' = Ê(m_1,
m_2,\dots,m_p)$, where the $m_i$ are positive integers for which $m_1 \leq m_2
\leq \dots \leq m_p$. ÊMoreover, if $m_{i-1} = m_i$ then $m_i < m_{i+1}$.

\smallskip

\item[ $ii)$ ] Given a pseudo type vector $T'$ and lines $L_1,\dots,L_p$, a
{\em pseudo linear configuration of type $T'$} is a set of points $\X = 
\X_1 \cup
\X_2 \cup\dots \cup \X_p$ where $\X_i$ is a set of $m_i$ points on 
$L_i$. Ê We do not allow any point of $\X_i$ to lie on line $L_j$ for $j 
\neq i$.

\smallskip

\item[iii)] A pseudo linear configuration of type $T' =(m_1,m_2,\ldots ,m_p)$
in $\PP^2$ is called Ê{\em standard} if it
consists of:
$$
\begin{array}{llcl}
m_p& \mbox{ points with coordinates }&[j:0:1]& 0\leq j\leq m_p -1, j \in
\mathbb N,\\
\\
ÊÊ Ê Ê & Ê Ê Ê& \vdots Ê Ê Ê Ê Ê Ê Ê Ê Ê Ê Ê Ê Ê Ê Ê Ê Ê Ê& Ê\\ \\
m_2& \mbox{ points with coordinates }&[j:p-2:1]& 0\leq j\leq m_2 -1, j \in
\mathbb N,\\ \\
m_1& \mbox{ points with coordinates }&[j:p-1:1]& 0\leq j\leq
m_1 -1, j \in \mathbb N.
\end{array}
$$

\end{itemize}
\end{definition}

We now describe an O-sequence that can be associated to a pseudo type vector
which depends only on the numerical information that is contained in the pseudo
type vector. We wish to stress, however, that we are not claiming that {\em
every} pseudo linear configuration with the given pseudo type vector has this
O-sequence as the first difference of its Hilbert function. ÊWe will see later
(see Theorem \ref{pseudo k con = bdl}) that such a strong statement is true
only for certain pseudo-type vectors.

\begin{definition} \label{standard os}
Let $T' = (m_1,\dots,m_p)$ be a pseudo type vector. ÊThe {\em standard
O-sequence } associated to $T'$ is given by a ``shifted sum'' of certain
sequences $s_i$ defined as follows: Êif we formally suppose that $m_0 = 0$ and
$m_{p+1} = \infty$, then

\begin{itemize}

\item If $m_{i-1} < m_i < m_{i+1}$ then
\[
(s_i)_t = \left \{
\begin{array}{ll}
1 & \hbox{for } 0 \leq t \leq m_i-1 ;\\
0 & \hbox{otheriwise}.
\end{array}
\right.
\]

\item If $m_{i-1} = m_i < m_{i+1}$ then
\[
(s_i)_t = \left \{
\begin{array}{ll}
1 & \hbox{for } t=0 ; \\
2 & \hbox{for } 1 \leq t \leq m_i-1 ; \\
1 & \hbox{for } t= m_i \\
0 & \hbox{otherwise}.
\end{array}
\right.
\]

\item If $m_{i-1} < m_i = m_{i+1}$ we do not define $s_i$.

\end{itemize}
Also, for a sequence $s_i$ and non-negative integer $k$, we define the shifted
sequence $s_i(-k)$ to be a rightward shift of $s_i$ by $k$ places (so instead
of starting in degree 0 it starts in degree $k$).

Then the {\it standard O-sequence associated to $T^\prime$} is:
\[
\sum_{i=1}^p s_i (i-p)
\]
where it is understood that the sum skips any indices for which $s_i$ is 
not defined.
\end{definition}

\begin{remark} \label{std hf from os}
The {\em standard Hilbert function associated to $T'$} is the numerical 
function whose first difference is the standard O-sequence associated to
$T'$ (as defined above). \qed
\end{remark}

\begin{example} \label{pseudo from 3 6 7}
Let $T' = (3,6,6,7,12,14)$. ÊThen the standard O-sequence associated to
$T^\prime$ is:

\medskip

\begin{center}

\begin{tabular}{cccccccccccccccc}
&&&&&1 & 1 & 1 \\
&&& 1 & 2 & 2 & 2 & 2 & 2 & 1 \\
&& 1 & 1 & 1 & 1 & 1 & 1 & 1 \\
& 1 & 1 & 1 & 1 & 1 & 1 & 1 & 1 & 1 & 1 & 1 & 1 \\
1 & 1 & 1 & 1 & 1 & 1 & 1 & 1 & 1 & 1 & 1 & 1 & 1 & 1 \\ \hline
1 & 2 & 3 & 4 & 5 & 6 & 6 & 6 & 5 & 3 & 2 & 2 & 2 & 1
\end{tabular}

\end{center}

\medskip

\noindent (Note that there is no $s_2$.) ÊThis sequence is the first 
difference
of the O-sequence:
\[
1 \ \ 3 \ \ 6 \ \ 10 \ \ 15 \ \ 21 \ \ 27 \ \ 33 \ \ 38 \ \ 41 \ \ 43 \ \ 45 \
\ 47 \
\ 48 \ \ 48 \dots
\]
which is the standard Hilbert function associated to $T'$. \qed
\end{example}

\begin{remark}\label{couldve}
We will see that Definition \ref{standard os} was designed so that the 
number of sequences $s_i$ correspond to the number of applications of basic 
double linkage. ÊThe sequences $s_i$ containing 2's will correspond to 
basic double links of the form $J = Q I + (F)$, where $Q = L_1 L_2$ is a 
product of linear forms.

It should be noted that such a basic double link can also be viewed as a
sequence of two basic double links $J_1 = L_1 I + (F)$ and $J = L_2 J_1 + (F)$
(with the same $F$). ÊBecause of this, we could also write the O-sequence sum
without any 2's. ÊFor instance, the above computation would become

\begin{center}

\begin{tabular}{cccccccccccccccc}
&&&&&1 & 1 & 1 \\
&&&& 1 & 1 & 1 & 1 & 1 & 1 Ê& \\
&&& 1 & 1 & 1 & 1 & 1 & 1 & Ê\\
&& 1 & 1 & 1 & 1 & 1 & 1 & 1 \\
& 1 & 1 & 1 & 1 & 1 & 1 & 1 & 1 & 1 & 1 & 1 & 1 \\
1 & 1 & 1 & 1 & 1 & 1 & 1 & 1 & 1 & 1 & 1 & 1 & 1 & 1 \\ \hline
1 & 2 & 3 & 4 & 5 & 6 & 6 & 6 & 5 & 3 & 2 & 2 & 2 & 1
\end{tabular}

\end{center}

\noindent and we would not have to worry about the extra shift. However, for
our purposes (constructing 2-fat point schemes) it is important to do the basic
double link in one step, so we retain this slightly more complicated 
notation. \qed
\end{remark}

We will see in Theorem \ref{pseudo k con = bdl} that there is a very 
precise condition that determines whether or not the Hilbert function of a 
pseudo linear configuration is uniquely determined by its 
type. ÊNevertheless, we now show that the Hilbert function of a {\em 
standard} pseudo linear configuration is uniquely determined, and in fact 
is equal to the function described in Definition \ref{standard os}.

\begin{proposition} \label{std pseudo has std os}
Let $\X \subset \mathbb P^2$ be a standard pseudo linear configuration of 
type $T' =$ \linebreak $ (m_1,m_2,\dots,m_p)$. ÊLet $\Delta T' = (m_1, 
m_2-m_1, \dots, m_p - m_{p-1})$. ÊThen:

\begin{itemize}
\item[{\it i)}] $\X$ can be built up by basic double linkage;

\item[{\it ii)}] the first difference of the Hilbert function of $\X$ is 
the standard O-sequence associated to $T'$;

\item[{\it iii)}] assume that between any two zero entries of $\Delta T'$ 
there is at least one entry $>1$. ÊIf $\Delta T'$ ends with a $0$, or Êwith 
a sequence
$(\dots,0,1,\dots,1)$ (i.e. a 0 followed by any number of 1's), then the 
regularity of
$\X$ is $m_p+1$. ÊOtherwise the regularity is $m_p$.

\item[{\it iv)}] If there are zero entries between which there are no 
entries $>1$ then the regularity of $\X$ may be arbitrarily larger than $m_p$.

\end{itemize}
\end{proposition}

\begin{proof}
Let $T' = (m_1,\dots,m_p)$, $T'' = (m_1,\dots,m_{p-1})$. ÊLet $\X$ be a 
standard
pseudo linear configuration of type $T'$ and let $\X_1$ be the obvious 
subset which is a pseudo linear configuration of type $T''$. ÊWe will show 
that $\X$ can be obtained from $\X_1$ by basic double linkage. ÊLet
\[
F = x(x-z)\dots(x-(m_{p-1}z)) \dots (x-(m_p -1)z).
\]
Then $F \in I_{\X_1}$ (since the configuration $\X$ is standard) and $\deg 
F = m_p$. ÊLet $G = y$.

Consider the basic double link
\[
I = G \cdot I_{\X_1} + (F).
\]
Then by Lemma \ref{BDL lemma}, $I$ is a saturated ideal whose support is 
exactly
$\X$. ÊTo show $I = I_\X$, it remains to show that $I$ is reduced. ÊBut the 
degree of the scheme defined by $I$ is $m_p$ more than $\deg \X_1$, by 
Lemma \ref{BDL lemma}, so $I$ and $I_\X$ are saturated ideals of 
zero-dimensional schemes with the same support and same degree, one of 
which is reduced. Hence they are equal. ÊThis proves $i)$.

By taking the first difference of (\ref{hf formula}), we obtain
\begin{equation} \label{what we need}
\Delta h_{R/I_\X} (t) = \Delta h_{R/(F,G)}(t) + \Delta h_{R/I_{\X_1}}(t-1).
\end{equation}
Using the fact that $\Delta
h_{R/(F,G)}(t)$ is
\[
\begin{array}{cccccccccccccc}
1 & 1 & \dots & 1 & 0 \\
(0) & (1) & \dots & (m_p -1) & (m_p)
\end{array}
\]
and taking Remark \ref{couldve} into account, it is clear that the first
difference of the Hilbert function of $\X$ is obtained (inductively) by 
\ref{standard os} and so $\X$ has the standard O-sequence associated to 
$T'$. ÊThis proves $ii)$.

We now verify the conclusions of $iii)$ by induction, assuming that it 
holds for
$\X_1$. ÊThe technical assumption, i.e. that there is at least one entry 
$>1$ between any two zero entries, remains true for $\X_1$. ÊNote that if 
$\Delta T''$ ends with a 0 or with a sequence $(\dots,0,1,\dots,1)$, then 
$\X_1$ has
regularity Ê$m_{p-1}+1$, by induction. ÊOtherwise $\X_1$ has regularity 
$m_{p-1}$.

\medskip\noindent {\bf Case 1.} Suppose that $\Delta T''$ ends with a 0 or 
with a sequence $(\dots,0,1,\dots,1)$, and that $m_p = m_{p-1}+1$. ÊThen 
the first difference of the Hilbert function of $\X_1$ ends in degree 
$m_{p-1}$ (because of the regularity noted above); but then (\ref{what we 
need}) forces this to be shifted by 1, giving precisely what is required in 
the O-sequence computation of Definition \ref{standard os}. ÊIn particular, 
the first difference of the Hilbert function of $\X$ ends in degree 
$m_{p-1}+1 = m_p$, and $\X$ has regularity
$m_p+1$ as claimed.

\medskip\noindent {\bf Case 2.} ÊSuppose that $\Delta T''$ ends with a 0 or 
with a sequence $(\dots,0,1,\dots,1)$, and that $m_p > m_{p-1}+1$. ÊThen 
again (\ref{what we need}) gives the O-sequence computation of Definition 
\ref{standard os}. Ê Again the first difference of the Hilbert function of 
$\X_1$, shifted by 1, Êends in degree $m_{p-1}+1 < m_p$, but Êthis time the 
regularity is determined by the $m_p$ new points, and is equal to $m_p$.

\medskip\noindent {\bf Case 3.} If $\Delta T''$ does not end with a 0 or 
with a sequence $(\dots,0,1,\dots,1)$, and if $m_p > m_{p-1}+1$ then the 
same argument as in Case 2 applies.

\medskip\noindent {\bf Case 4.} If $m_{p-1} = m_p$ then necessarily we have 
$m_{p-2}<m_{p-1}$, by the definition of a pseudo linear 
configuration. ÊThus, for the pseudo type vector $T''$, $s_{p-1}$ is 
defined. ÊHowever, for $T'$, since $m_{p-2} < m_{p-1} = m_p$, we obtain 
that $s_{p-1}$ is not defined, but $s_p$ is of the second type described in 
the O-sequence computation of Definition \ref{standard os}. ÊThus, in this 
case, the induction takes us from $T''' = (m_1,\dots,m_{p-2})$ to $T' = 
(m_1,\dots,m_p)$. ÊBut now the technical assumption that between any two 
zero entries of $\Delta T'$ there is at least one entry $>1$, together with 
induction, guarantees that the regularity of the standard pseudo linear 
configuration $\X_2$
determined by $T'''$ is $\leq m_p-1$. ÊThe O-sequence computation of 
Definition \ref{standard os} indicates that we must shift the O-sequence of 
$\X_2$ by 2, so that it now ends in degree $\leq m_p+1$. ÊHence the
regularity is computed by the two sets of $m_p$ points, and is equal to
$m_p+1$.

For $iv)$ it is enough to realize that in the standard configuration of type
\[
(1,1,2,2,3,3,\dots,m,m)
\]
there is a set of $2m$ points lying on the ``vertical'' line $x$, so the 
regularity is at least $2m$.
\end{proof}

The standard pseudo linear configuration is clearly very 
special. ÊNevertheless, we now show that the technical assumption that 
between any two zero entries there is at least one entry $>1$ (used to 
control the regularity of the standard pseudo linear configuration), is 
enough to guarantee that {\em all} pseudo linear configurations of that 
type have the same Hilbert function.

\begin{theorem} \label{pseudo k con = bdl}
Consider a pseudo type vector $T' = (m_1,m_2,\dots,m_p)$. ÊLet
\[
\Delta T' = (m_1-0, m_2-m_1, \dots, m_p - m_{p-1})
\]
be its first difference (note that $\Delta T'$ has all entries non-negative).
Then {\bf every} pseudo linear configuration of type $T'$ can be realized 
as the
result of a sequence of basic double links if and only if the following
condition holds:

\begin{equation} \label{condition}
\begin{hbox}
{ Between any two zero entries of $\Delta T'$ there is at least one entry that
is $>1$.}
\end{hbox}
\end{equation}

If this condition holds then the Hilbert function of any pseudo linear 
configuration of type $T^\prime$ is the same. ÊThe first difference of that 
Hilbert function is the O-sequence given by Definition \ref{standard os}.

In particular, if condition (\ref{condition}) holds for the pseudo type 
vector, $T'$, of a given pseudo linear configuration, then the regularity 
of that pseudo linear configuration is determined as follows: if $\Delta 
T'$ ends with a 0, or with a sequence $(\dots 0,1,1,\dots,1)$ (i.e.\ a 0 
followed by any number of 1's), then the regularity is $m_p +1$. ÊOtherwise 
the regularity is $m_p$.

If condition (\ref{condition}) does not hold then the Hilbert function of a
pseudo linear configuration of type $T^\prime$ is not uniquely determined.
\end{theorem}

\begin{proof}
Note that by the definition of a pseudo type vector, there must be at least one
non-zero entry between any two zeroes in the vector $\Delta T'$. Ê We first
prove that the condition (\ref{condition}) is sufficient Êto realize a given
pseudo linear configuration as being obtained as a sequence of basic double
linkages, by working from left to right in the pseudo type vector (imagine 
a pointer moving along the marker and keeping track of our current position).

Having built the subconfiguration corresponding to entries 
$m_1,\dots,m_{i-1}$,
the next step:

\begin{itemize}

\item will involve only $m_i$ if $m_i < m_{i+1}$,

\item Êwill involve $m_i$ and $m_{i+1}$ if $m_i = m_{i+1}$.

\end{itemize}

We will use the fact that sets consisting of $m_i$ points on a line, or
$m_i$ points on each of two lines (avoiding the intersection point of the 
two lines) are both complete intersections in $\P^2$. Ê(This is no longer 
necessarily true for three lines.)

Of course, if $T'$ satisfies (\ref{condition}) then so does every 
subsequence.
To simplify the notation, at each step we will take: $\X$ to be the
subconfiguration built up so far (by induction); $\Y$ to be the set added; and
$\Z$ will be the new set, $\Z = \X \cup \Y$. ÊNote that if $m_i < m_{i+1}$ 
then $\Y$ is a set of $m_i$ points on a line $L$ (and by abuse of notation 
we denote by $L$ also the linear form defining this line), and if $m_i = 
m_{i+1}$ then $\Y$ consists of $m_i$ points on each of two lines, and we 
denote by $Q$ this union of lines (and the corresponding product of two 
linear forms).

To begin the construction we take $\X$ to be a set of

\begin{itemize}

\item $m_1$ points on line $L_1$ if $m_1 < m_2$,

\item $m_1$ points on each of lines $L_1$ and $L_2$ if $m_1 = m_2$, 
avoiding the intersection point of $L_1$ and $L_2$.

\end{itemize}

In the first case $I_\X$ has generators of degrees $1$ and $m_1$, and 
regularity
$m_1$. In the second case $I_\X$ has generators of degrees 2 and $m_1$, and
regularity $m_1 +1$.

\medskip Now let $\X$ be the configuration constructed up to entry 
$m_{i-1}$. ÊNote that we necessarily have $m_{i-1} < m_i$, since if they 
were equal then we would have constructed the points corresponding to 
$m_{i-1}$ and $m_i$ at the same time. ÊWe have the partial first difference 
vector $(m_1 - 0, m_2-m_1, \dots, m_{i-1}-m_{i-2})$. ÊBy induction, if this 
first difference vector ends with a 0 or with a sequence $(\dots, 
0,1,1,\dots,1)$ then $\reg(I_\X) = m_{i-1}+1$, and otherwise $\reg(I_\X) 
= Êm_{i-1}$.

\medskip\noindent {\bf Case 1:} $m_i < m_{i+1}$.

This means that we want to add $m_i$ points on $L$. ÊWe have an exact sequence
\begin{equation} \label{std ex seq}
\begin{array}{cccccccccccccccccc}
0 & \rightarrow & [I_\Z : L](-1) & \stackrel{\times L}{\longrightarrow} &
I_\Z & \rightarrow & \displaystyle \frac{I_\Z + (L)}{(L)} & \rightarrow & 0. \\
&& || \\
&& I_\X (-1)
\end{array}
\end{equation}
We sheafify and take cohomology. ÊNote that $\displaystyle \widetilde{\left (
\frac{I_\Z + (L)}{(L)} \right )} = {\mathcal I}_{\Y|L}$ is the ideal sheaf 
of $\Y$, viewed as a subscheme of $L = \mathbb P^1$. ÊIts global sections 
begin in degree $m_i$. ÊWe know that the regularity of $I_\X$ is $m_{i-1}$ 
or $m_{i-1}+1$, so $h^1 (\mathcal I_\X (m_{i-1})) = 0$. ÊNote that $m_{i-1} 
\leq m_i-1$. ÊCondition (\ref{condition}) does not directly affect this 
case since we have assumed that it holds for $\X$ and we have $m_{i-1} < m_i$.

ÊFrom the exact sequence
\[
0 \rightarrow \mathcal O_{\mathbb P^2}(-1) \rightarrow \mathcal I_\Y 
\rightarrow
\mathcal I_{\Y|L} \rightarrow 0
\]
we get $h^1(\mathcal I_\Y(t)) = h^1(\mathcal I_{\Y|L}(t))$ for all $t \geq 
-1$.
Sheafifying (\ref{std ex seq}), twisting by $t \geq 0$ and taking 
cohomology, we
get
\begin{equation} \label{std in cohom}
0 \rightarrow (I_\X)_{t-1} \rightarrow (I_\Z)_t \rightarrow (I_{\Y|L})_t 
\rightarrow
H^1(\mathcal I_\X(t-1)) \rightarrow H^1(\mathcal I_\Z(t)) \rightarrow 
H^1(\mathcal
I_\Y(t)) \rightarrow \dots
\end{equation}
Taking $t = m_i$, we have $H^1(\mathcal I_\X(m_i -1)) = 0$ since $m_{i-1} 
\leq m_i -1$ and $h^1(\mathcal I_\X(m_{i-1})) = 0$. ÊHence the restriction 
$(I_\Z)_{m_i} \rightarrow (I_{\Y|L})_{m_i}$ is a surjection, and the 
non-zero element of $(I_{\Y|L})_{m_i}$ lifts to a form $F \in I_\Z$ that 
does not vanish on the line $L$, so in particular (since the points of $\Y$ 
are distinct) meets
$L$ transversally in $\Y$. ÊSince $\Z = \X \cup \Y$ and $\Y$ is the complete
intersection of $F$ and $L$, $\Z$ is a basic double link: indeed, $I_\Z$ and $L
\cdot I_\X + (F)$ are saturated ideals defining the same set of points, so 
we have $I_\Z = L \cdot I_\X + (F)$.

We now verify the Hilbert function calculation. ÊSince we know that $I_\Z = L
\cdot I_\X + (F)$, we can use the theory of basic double linkage as 
described in
Section \ref{prelim section}. ÊIndeed, it follows easily from (\ref{hf 
formula}). ÊLet $G = L$, Ê$d_1 = m_i$ and $d_2 = 1 \ (= \deg L)$. ÊWe then 
note that we are in the first case of the O-sequence computation of 
Definition \ref{standard os}, and that in that formula $(s_i)$ is now just 
the first difference of $h_{R/(F,L)}(t)$. ÊThen the bottom row of the 
O-sequence computed in Definition \ref{standard os} (see Example 
\ref{pseudo from 3 6 7}) corresponds to the first difference of the Hilbert 
function of $R/(F,L)$, and the rows above
the bottom row all together correspond to (a decomposition of) the first
difference of the Hilbert function of $\X$, shifted by 1. ÊThe connection 
is made
by (\ref{hf formula}).

ÊÊAs for the regularity, we know that
\[
\reg (I_\Z) = 1 + \min \{ t \ | \ h^1 (\mathcal I_\Z(t)) = 0 \}.
\]
The sequence (\ref{std in cohom}) shows that $\reg (I_\Z)$ is the larger 
of Ê$\reg (I_\X)+1$ Êand $m_i$, depending (respectively) on whether $m_i = 
m_{i-1} +1$ or $m_i > m_{i-1}+1$. ÊUsing induction, this shows that 
$\reg(I_\Z) = m_i + 1$ if the first difference vector $(m_1-0, m_2 -m_1, 
\dots , m_i - m_{i-1})$ ends with a sequence $(\dots 0,1,\dots,1)$ (in this 
Case, it is excluded that this first difference will end with a 0) Êand 
$\reg (I_\X) = m_i$ otherwise.

\medskip\noindent {\bf Case 2.} $m_{i} = m_{i+1} $.

In this case we let $\X$ be the set of points corresponding to the pseudo 
type vector $(m_1,\dots,m_{i-1})$, and we add $m_i$ points on each of two 
lines. ÊIn our argument, instead of $L$ we use $Q$, the Êunion of the two 
lines containing $m_i$
points each. ÊIn this Case the first difference vector $(m_1-0, m_2-m_1, \dots,
m_{i+1}-m_i)$ ends with a 0, so condition (\ref{condition}) implies that it 
does not end $(\dots, 0,1,\dots,1,0)$ (with only 1's between the 0's). ÊThe 
two possibilities are that {\bf (a)} there is no other 0, or else {\bf (b)} 
there is at least one entry that is $>1$ between the 0's.

We first would like to compute the regularity of $I_\X$. ÊIf {\bf (a)} 
holds then $m_1 < m_2 < \dots < m_{i-1} < m_i = m_{i+1}$. ÊHence $\X$ is a 
linear configuration, and its regularity is $m_{i-1}$. ÊWe remark that in 
this case $\Z$ is in fact a linear configuration minus a point on the 
``longest'' row, so its Hilbert function is just the truncation of the 
Hilbert function of the linear configuration of type 
$(m_1,\dots,m_{i-1},m_i, m_i+1)$ (cf.\ \cite{sindipaper}).

If {\bf (b)} holds, then there are again two possibilities. ÊFirst, it 
could happen that the first difference vector $(m_1-0,\dots,m_{i-1} - 
m_{i-2})$ ends in a 0. ÊIn this case $m_{i-1} = m_{i-2}$, and by induction 
the regularity of $I_\X$ is $m_{i-1} +1$. ÊHowever, condition {\bf (b)} 
then means that $m_i \geq m_{i-1} +2$.

The other possibility in {\bf (b)} is that the first difference vector
$(m_1-0,\dots,m_{i-1} - m_{i-2})$ does not end in a 0. ÊIf it ends in a string
$(\dots, 0, 1, \dots, 1)$ (all 1's after the 0) then again the regularity 
of $I_\X$ is $m_{i-1} +1$ by induction, but again {\bf (b)} forces $m_i 
\geq m_{i-1} +2$. ÊIf the first difference vector does not end in such a 
string, then by induction the regularity of $I_\X$ is $m_{i-1}$.

We conclude from the above analysis that in every case,
\begin{equation} \label{condition on h1}
h^1({\mathcal I}_\X (m_i-2)) = 0.
\end{equation}

By analogy with Case 1, but now using $Q$ instead of $L$, the exact
sequence (\ref{std ex seq}) now becomes
\[
\begin{array}{cccccccccccccccccc}
0 & \rightarrow & [I_\Z : Q](-2) & \stackrel{\times Q}{\longrightarrow} &
I_\Z & \rightarrow & \displaystyle \frac{I_\Z + (Q)}{(Q)} & \rightarrow & 0. \\
&& || \\
&& I_\X (-2)
\end{array}
\]
Sheafifying, twisting by $m_i$ and taking cohomology, we get
\[
0 \rightarrow (I_\X)_{m_i -2} \rightarrow (I_\Z)_{m_i} \rightarrow 
(I_{\Y|Q})_{m_i}
\rightarrow H^1(\mathcal I_\X(m_i -2)) \rightarrow H^1(\mathcal I_\Z(m_i)) 
\rightarrow
ÊÊ\dots
\]
By (\ref{condition on h1}), we have $h^1(\mathcal I_\X(m_i -2)) =0$. ÊSince 
$\Y$
consists of $m_i$ points on each of the two components of $Q$, and since by 
Definition \ref{def of pseudo} we do not allow a point of $\Y$ to lie on 
both components of $Q$, we first claim that a non-zero element of 
$(I_{\Y|Q})_{m_i}$ cannot vanish identically on either component of 
$Q$. ÊIndeed, if it vanished on either component then it would lift to a 
homogeneous polynomial $F \in (I_\Z)_{m_i}$ vanishing on a line, which then 
has as a factor a form $G$ of degree $m_i-1$ that vanishes on $m_i$ points 
on the other component of $Q$, but does not vanish identically on that 
component (since it is a non-zero element of $I_{\Y|Q}$). ÊImpossible.

Thus the vanishing of the first cohomology and the fact that $I_{\Y|Q}$ 
begins in degree $m_i$ (recall that $\Y$ is a complete intersection of type 
$(2, m_i)$) means that there is a form $F \in I_\Z$ of degree $m_i$ that 
does not vanish on either component of $Q$, and cuts out $\Y$ on $Q$ (in 
particular it meets $Q$ transversally). ÊSince $\Z = \X \cup \Y$, we again 
recognize $\Z$ as being obtained as a basic double link from $\X$, and we have
\[
I_\Z = Q \cdot I_\X + (F).
\]

By an argument similar to the one above, we can compute the Hilbert 
function of $\Z$ using the O-sequence computation in Definition 
\ref{standard os}. ÊIn this case the bottom row is given by the second case 
in that computation, since $\Y$
is a complete intersection of type $(2, m_i)$, and the shift between the 
bottom row of the computation and the rows above it is now 2 (see Example 
\ref{pseudo from 3 6 7}).

In fact, starting from (\ref{bdl exact seq}) we can easily compute the 
minimal free resolution of $I_\Z$, using a mapping cone and using a minimal 
free resolution
\[
0 \rightarrow \mathbb F_2 \rightarrow \mathbb F_1 \rightarrow I_\X 
\rightarrow 0.
\]
We get a free resolution
\[
0 \rightarrow
\begin{array}{c}
\mathbb F_2 (-2)\\
\oplus \\
R(-m_i-2)
\end{array}
\rightarrow
\begin{array}{c}
\mathbb F_1 (-2) \\
\oplus \\
R(-m_i)
\end{array}
\rightarrow I_\Z \rightarrow 0.
\]
However, Êby (\ref{condition on h1}) we have that $\reg(I_\X) \leq m_i-1$, 
so $F$
(having degree $m_i$) cannot be a minimal generator of $I_\X$. ÊBut the 
resolution is minimal if and only if $F$ is not a minimal generator (Lemma 
\ref{BDL lemma} (e)). Ê In particular, it follows that $\reg(I_\Z) = m_i 
+1$. Ê(Note that in this Case the sequence $\Delta T'$ ends with a 0, so 
this is the regularity claimed in the statement of the theorem.)

This completes one direction of the theorem. ÊFor the converse, we have to 
show that if $\Delta T'$ contains a subsequence $(\dots, 0, 1, \dots, 1,0)$ 
(all 1's between the two 0's) then there exist (at least) two pseudo linear 
configurations of type $T'$ with different Hilbert functions. ÊTo do this, 
first we will show that if $\Delta T'$ {\em ends} with such a subsequence, 
with no such subsequence preceeding it, then the conclusion holds. ÊSecond, 
we will show the general statement. ÊNote that we showed in Proposition 
\ref{std pseudo has std os} that for {\em any} pseudo type vector, there 
always exists one pseudo linear configuration (the standard one) that can 
be constructed by basic double linkage, and hence has ÊHilbert function 
whose first difference is given by the O-sequence computed in Definition 
\ref{standard os}. ÊSo for both the first part and the second part, we have 
to show that such a subsequence allows for a pseudo linear configuration 
that can {\em not} be constructed entirely by basic double links.

Suppose that $\Delta T'$ ends with the subsequence $(\dots, 0,1,\dots,1,0)$ 
(all 1's between the 0's) and no such subsequence precedes it. ÊThis means 
that if $T' = (m_1,\dots,m_{p-2}, m_{p-1},m_p)$ then $m_p = m_{p-1} = 
m_{p-2}+1$. ÊIn this paper we usually handle the case where $m_{p-1} = m_p$ 
by doing a basic double link using a quadric form $Q$, as described above 
(because of the application to non-reduced schemes that will be given 
below). However, just for this step in the current proof, it is convenient 
to view it as two separate basic double links using linear forms.

Let $T''$ be the pseudo type vector $(m_1,\dots, m_{p-2}, m_{p-1})$. ÊThen 
$\Delta T'$ ends in a sequence $(\dots, 0, 1, \dots, 1)$ where the end 
consists of nothing but 1's. Ê We have assumed that $T''$ satisfies 
(\ref{condition}). ÊLet $\X$ be a pseudo linear configuration corresponding 
to $T''$. ÊIt can be constructed by basic double linkage, and its Hilbert 
function is as described. ÊFurthermore, it follows from what we have 
already proven that the regularity of $I_\X$ is $m_{p-1}+1 = Êm_{p}+1$.

Now consider an additional line $L_p$, and choose $\Y$ to be a general set 
of $m_p$ points on $L_p$. ÊLet $\Z = \X \cup \Y$. Ê $\Z$ is a basic double 
link of $\X$ if and only if there is a form $F \in (I_\X)_{m_p}$ that 
contains $\Y$ but does not vanish on $L_p$.

We have an exact sequence
\[
0 \rightarrow {\mathcal I}_\X (m_p-1) \stackrel{\times 
L_p}{\longrightarrow} {\mathcal
I}_\X(m_p)
\rightarrow {\mathcal O}_{L_p} (m_p) \rightarrow 0.
\]
Since the regularity of $I_\X$ is $m_p+1$, we have
\begin{equation} \label{eexx}
0 \rightarrow (I_\X)_{m_p-1} \rightarrow (I_\X)_{m_p} 
\stackrel{r}{\longrightarrow}
H^0({\mathcal O}_{L_p} (m_p)) \rightarrow H^1({\mathcal I}_\X(m_p-1)) 
\rightarrow 0
\end{equation}
where the last cohomology group is not zero. ÊChoosing $\Y$ as above is 
equivalent to choosing a general element of the vector space $H^0({\mathcal 
O}_{L_p}(m_p))$. ÊThe image of $r$ is a proper subspace of $H^0({\mathcal 
O}_{L_p}(m_p))$, so the general section of $H^0({\mathcal O}_{L_p}(m_p))$ 
defining $\Y$ is not in the image of $r$. ÊWe conclude that any form in
$(I_\X)_{m_p}$ that vanishes on $\Y$ must in fact vanish on all of 
$L_p$. ÊHence we cannot express $\Z$ as a basic double link of $\X$.

We claim that the value of the Hilbert function of such a $\Z$ differs, in 
degree $m_p$, from the value of the corresponding Hilbert function given by the
O-sequence computation in Definition \ref{standard os} (see Remark \ref{std 
hf from os}). ÊIndeed, suppose that $\Z'$ were a pseudo linear 
configuration of the same type that was produced by basic double linkage, and
hence has the standard Hilbert function for that type. Ê(It is not hard to 
check that we can even assume that $\Z'$ is built up from the same $\X$, 
choosing the points of $\Y$ in a more careful way.) ÊIn degree $m_p$, the 
forms that vanish on $\Z$ consist entirely of products of $L_p$ with forms 
of degree $m_p-1$ vanishing on $\X$ (as discussed above), while $\Z'$ has 
those but also has a form of degree $m_p$ that is not of that form. ÊHence 
the Hilbert functions differ in degree $m_p$.

Now we prove the second part. ÊLet $T'$ be a pseudo linear configuration 
not satisfying (\ref{condition}). ÊIts first difference has an initial 
subsequence with first difference $(\dots, 0,1,1,\dots,1,0)$, where the 
earlier entries do satisfy (\ref{condition}). ÊLet $\Z$ and $\Z'$ be as 
above, both pseudo linear configurations with type given by this 
subsequence, and having different Hilbert functions. ÊWe claim that term by 
term we can add $m_i$'s to the subsequence, and correspondingly add points 
on a line that arise by basic double linkage. Ê(We do not claim that {\em 
only} basic double links are possible if there is another subsequence 
$(\dots, 0,1,1,\dots,1,0,\dots)$, but
only that in particular a basic double link {\em is} possible.) ÊThe basic 
double
links at each step are numerically the same, so they add the same amount in 
each
degree to the Hilbert functions. ÊSince we started with $\Z$ and $\Z'$ 
having different Hilbert functions, this will say that at each step the 
results have different Hilbert functions, and we will be finished.

We have seen that for any type there exists a standard pseudo linear 
configuration, so we can assume that $\Z'$ is a standard pseudo linear 
configuration, and it can continue to be built up by basic double links as 
claimed. ÊThe real assertion is that this is true of $\Z$ as well. ÊFirst 
note that if we twist the exact sequence (\ref{eexx}) by any $t>0$, we 
obtain the short exact sequence
\[
0 \rightarrow (I_\X)_{m_p+t-1} \rightarrow (I_\X)_{m_p+t} 
\stackrel{r}{\longrightarrow}
H^0({\mathcal O}_{L_p} (m_p+t)) \rightarrow 0
\]
because of the regularity. Now choose $t$ so that $m_p+t = m_{p+1}$ (recall 
that both $\Z$ and $\Z^\prime$ ended with $m_p = m_{p-1}$). ÊThis says, in 
particular, that there is a form $F$ in $(I_\X)_{m_{p+1}}$ vanishing on 
$\Z$ but not vanishing identically along $L_p$, because $r$ is surjective 
and we can choose a section of $H^0({\mathcal O}_{L_p}(m_{p+1}))$ that 
vanishes at $\Y$ plus $t$ general points of $L_p$. ÊBut then choosing a 
general line $L_{p+1}$, this meets the same $F$ in $m_{p+1}$ distinct 
points, forming a basic double link of $\Z$ of type 
$(m_1,\dots,m_p,m_{p+1})$. ÊNow it is trivial to build up the rest of the 
pseudo type by basic double linkage, since we can take products of $F$ with 
general forms of suitable degree to produce the points.
\end{proof}

\begin{example}
Theorem \ref{pseudo k con = bdl} applies to the pseudo type vectors
\[
(2,2, 3,4,5,7,7), \ (2,2,3,4,6,7,7), \ (2,2,3,4,5,6,7), \ (5,5,7,7)
\]
but not to
\[
(2, 2, 3, 4, 5, 6, 7, 7) \hbox{ or } (6,6,7,7).
\]
A simpler example to show that the Hilbert function may vary if 
(\ref{condition}) does not hold is the pseudo type vector $(1,1,2,2)$. ÊIf 
these points form a standard pseudo linear configuration, i.e.\
\[
\begin{array}{ccccccccccc}
\bullet \\
\bullet \\
\bullet & \bullet \\
\bullet & \bullet
\end{array}
\]
then the Hilbert function of the points is $(1,3,5,6,6,\dots)$ (note that 
by considering the ``vertical'' lines, this set of points is realized as a 
linear configuration of type $(2,4)$). ÊOn the other hand, if the points 
are chosen generically on the four ``horizontal'' lines then they are in 
fact 6 generic points in $\mathbb P^2$, so the Hilbert function is 
$(1,3,6,6,\dots)$. \qed
\end{example}


\section{The resolution of the ideal of a pseudo linear configuration}

In the last section we saw the necessary and sufficient condition for the 
Hilbert
function of a pseudo linear configuration to be uniquely determined from 
the type. ÊThis was seen to be equivalent to the condition that every 
pseudo linear configuration of given type can be built up by basic double 
linkage in the way prescribed by the type. ÊThis is analogous to the 
situation for linear configurations, where the type uniquely determines the 
Hilbert function (but with no condition needed).

For linear configurations, in fact, the type uniquely determines the graded 
Betti
numbers (which are maximal among all algebras with the given Hilbert function
\cite{GHSAdv}). ÊWe now turn to the question of when the type of a pseudo
linear configuration uniquely determines the graded Betti numbers, and how 
to determine what those graded Betti numbers are. ÊWe will use the fact 
that in $\mathbb P^2$, when the Hilbert function is fixed, the graded Betti 
numbers depend only on the degrees of the minimal generators.

\begin{example} \label{not unique}
Consider the pseudo type vector $(1,2,2,3)$. ÊWe have seen that any pseudo
linear configuration of this type arises from a sequence of two basic 
double linkages (starting from a single point), so the Hilbert function is 
uniquely determined. ÊBut we will see now that the graded Betti numbers are 
not uniquely determined. ÊIn particular, we will see that the form $F$ of 
degree 3 that is used for the last basic double linkage may or may not be a 
minimal generator of the subconfiguration of type $(1,2,2)$.

First suppose that the pseudo linear configuration is standard:
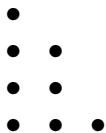
\[
\begin{array}{cccccccccccc}
\bullet \\
\bullet & \bullet \\
\bullet & \bullet \\
\bullet & \bullet & \bullet
\end{array}
\]
$\X$ consists of the five points on the first three ``horizontal'' 
lines. ÊIn this case $F$ can be taken to be the product of the three 
``vertical'' lines. ÊNote that the product of the leftmost two vertical 
lines is an element of the ideal of $\X$ (in fact it is the only generator 
of $I_\X$ of degree 2), so $F$ is not a minimal generator of $I_\X$. ÊHence 
no splitting occurs, by Lemma \ref{BDL lemma} (e). ÊWe have a minimal free 
resolution
\[
0 \rightarrow \begin{array}{c}
R(-5)^2 \\
\oplus \\
R(-4)
\end{array}
\rightarrow
\begin{array}{c}
R(-3) \oplus R(-4)^2 \\
\oplus \\
R(-3)
\end{array}
\rightarrow I_\Z \rightarrow 0.
\]

On the other hand, suppose that our pseudo linear configuration of pseudo type
$(1,2,2,3)$ is formed by general points on each of the four lines. ÊNote 
that $I_\X$ still has only one quadric generator, and this quadric meets 
the fourth
``horizontal'' line in two points, say $P_1$ and $P_2$. ÊSince $\Z$ was 
chosen with general points on each of the lines, the three points on this 
fourth ``horizontal'' line are disjoint from $P_1$ and $P_2$. ÊTherefore no 
$F$, cutting out the three points on this line, is Êa multiple of the 
quadric generator of $I_\X$; hence any such $F$ can be chosen as a minimal 
generator of degree 3. ÊTherefore a copy of $R(-4)$ splits off in the above 
resolution, for this pseudo linear configuration, and we obtain the minimal 
free resolution
\[
0 \rightarrow R(-5)^2 \rightarrow R(-3)^2 \oplus R(-4) \rightarrow I_\Z 
\rightarrow 0.
\]
Therefore, as claimed, the graded Betti numbers are not uniquely determined 
for the pseudo type vector $(1,2,2,3)$. \qed
\end{example}

With this example in mind, we consider the graded Betti numbers of a pseudo
linear configuration that arises as a result of basic double linkage (e.g.\ 
by satisfying condition (\ref{condition}) or by being a standard pseudo 
linear configuration). Ê Suppose that $I_\X$ has minimal free resolution
\[
0 \rightarrow \mathbb F_2 \rightarrow \mathbb F_1 \rightarrow I_\X 
\rightarrow 0
\]
and that $\Z$ arises from $\X$ by basic double linkage using $F$ and $L$ as 
before. ÊThen the diagram
\[
\begin{array}{cccccccccccccccccc}
&&&&& 0 \\
&&&&& \downarrow \\
&& 0 &&\mathbb F_2 (-1) & \oplus & 0\\
&& \downarrow &&& \downarrow \\
&& R(-m_i-1) && \mathbb F_1 (-1) & \oplus & R(-m_i) \\
&& \downarrow &&& \downarrow \\
0 & \rightarrow & R(-m_i-1) & \rightarrow & I_\X (-1) & \oplus & R(-m_i) & 
\rightarrow
& I_\Z & \rightarrow & 0
\end{array}
\]
ÊÊyields a resolution (using the mapping cone)
\[
0 \rightarrow
\begin{array}{c}
\mathbb F_2 (-1)\\
\oplus \\
R(-m_i-1)
\end{array}
\rightarrow
\begin{array}{c}
\mathbb F_1 (-1) \\
\oplus \\
R(-m_i)
\end{array}
\rightarrow I_\Z \rightarrow 0.
\]
As mentioned before, this resolution is minimal if and only if $F$ is {\em 
not} a
minimal generator of
$I_\X$ (Lemma \ref{BDL lemma} (e))

So we are reduced to the problem of determining whether or not $F$ is a minimal
generator of $I_\X$. ÊIf $m_i > \reg(I_\X)$ then clearly $F$ is not a 
minimal generator of $I_\X$. ÊIf $m_i = \reg(I_\X)$ (the only other 
possibility) then we know that $m_i = m_{i-1} +1$, and that the first 
difference of the partial type vector $(m_1,m_2,\dots,m_{i-1})$ ends with 
either a 0 or a sequence $(0,1,\dots,1)$. ÊHowever, unfortunately in this 
case the question of whether or not $F$ is a minimal generator of $I_\X$ is 
not merely a numerical one.

\begin{example} \label{122445}
One can check that the pseudo type vector $(1,2,2,4,4,5)$ does not have the 
property that all pseudo linear configurations of this type have the same 
graded Betti numbers. ÊIndeed, letting $\X_1$ be a point, with ideal 
$(A_1,A_2)$ ($\deg A_i = 1$), then we successively form
\[
\begin{array}{rcll}
I_{\X_2} & = & (QA_1, QA_2, F_1) & \hbox{where $F_1 \in I_{\X_1},
\deg F_1 = 2$, $\deg Q = 2$} \\
I_{\X_3} & = & (Q'QA_1, Q'QA_2, Q'F_1,F_2) & \hbox{where $F_2 \in 
I_{\X_2}$, $\deg F_2
= 4$,
$\deg Q' = 2$} \\
\end{array}
\]

The point is that $I_{\X_3}$ does have generators of degree 5, so forming the
last basic double link using a form $F \in I_{\X_3}$ of degree 5 can be 
done with $F$ a minimal generator of $I_{\X_3}$ or not. \qed
\end{example}

\medskip Theorem \ref{pseudo k con = bdl} gave (in particular) a necessary 
and sufficient condition for the Hilbert function of a pseudo linear 
configuration, $\X$, to be uniquely determined by the pseudo type; namely, 
(\ref{condition}), that between any two zero entries of $\Delta T'$ there 
is at least one entry that is $>1$. ÊWe would like to do the same thing for 
the graded Betti numbers. ÊOf course we have to begin by assuming 
(\ref{condition}), since if the Hilbert function can vary then so can the 
graded Betti numbers. ÊIn particular, we can assume that $\X$ can be 
realized as a sequence of basic double links.

The following lemma is trivial, but we will refer to it several times in the
next result.

\begin{lemma} \label{gen lemma}
Let $\Z$ be a basic double link of $\X$, so that $I_\Z = A \cdot I_\X + 
(F)$ with $F \in I_\X$. ÊAssume that $F$ is not a minimal generator of 
$I_\X = (G_1,\dots,G_r)$. ÊAssume further that the maximal degree of a 
minimal generator of $I_\X$ is $d$.

\begin{itemize}
\item[(a)] If $\deg A = 1$ then the minimal generators of $I_\Z$ have 
degrees Ê$\deg G_1 +1,\dots,$ \linebreak $\deg G_r+1, \deg F$. ÊIf $d \leq 
\deg F -1$ then all generators have degree $\leq \deg F$.

\medskip

\item[(b)] If $\deg A = 2$ then the minimal generators of $I_\Z$ have 
degrees $\deg G_1 +2,\dots,$ \linebreak $ \deg G_r +2, \deg F$. ÊIf $d \leq 
\deg F -2$ then all generators have degree $\leq \deg F$.
\end{itemize}

\end{lemma}

\begin{remark}
In the following theorem, we will be constructing a pseudo linear 
configuration $\Z$ inductively from a smaller one, $\X$, and studying the 
question of whether the polynomial $F$ used in the basic double link is a 
minimal generator of $I_\X$ or not. ÊIn each case, $F$ will have the 
largest possible degree allowed by the regularity. ÊIf our analysis shows 
that $I_\X$ does have a minimal generator of that degree, then a general 
element of $I_\X$ of that degree can form part of a minimal generating 
set. ÊTherefore, even though the argument that we use to show that $I_\X$ 
has a minimal generator of that degree will produce $F$ having components 
in common with other generators, it is just the existence of generators 
that is important, and then a general choice will have no such common 
components. \qed
\end{remark}

In the following result, it is helpful to keep in mind Example \ref{not 
unique} and Example \ref{122445}.

\begin{theorem} \label{betti unique}
Consider a pseudo type vector $T' = (m_1,m_2,\dots,m_p)$. ÊLet
\[
\Delta T' = (m_1-0, m_2-m_1, \dots, m_p - m_{p-1})
\]
ÊÊbe its first difference, and assume that (\ref{condition}) holds, i.e.\ 
between any two zero entries of $\Delta T'$ there is at least one entry 
that is $>1$. Let Ê$\Z$ be a pseudo linear configuration of pseudo type 
$T'$. ÊThen the following hold.

\begin{itemize}

\item[(a)] The graded Betti numbers of Ê$I_\Z$ are uniquely determined if 
and only if Ê$\Delta T'$ contains none of the following as subsequences:
\begin{equation}\label{bad list}
\begin{array}{c}
(1,0,1), \\
(1,0,2,0,1) \\
(1,0,2,0,2,0,1),\\
\vdots \\
(1,0,2,0,\dots,0,2,0,1)
\end{array}
\end{equation}

\item[(b)] If $\Delta T'$ contains none of (\ref{bad list}) as subsequences 
then the number of minimal generators of $I_\Z$ is $p+1-a$ where $a$ is the 
number of 0's appearing in $\Delta T'$.

\item[(c)] In particular, $\Z$ has the maximum number of minimal generators 
allowed by the Hilbert function if and only if it is a linear configuration 
(i.e. $\Delta T'$ contains no 0's).

\end{itemize}
\end{theorem}

\begin{proof}
We know that $\Z$ can be obtained by a sequence of basic double links, since
(\ref{condition}) holds. ÊAt each step the ideal has the form $J = AI + 
(F)$ where $A$ is a form of degree 1 or 2, $F \in I$, and $(A,F)$ is a 
regular sequence. ÊIf $I = (F_1,\dots,F_r)$ then $J$ is generated by $(F, 
AF_1, \dots, AF_r)$. ÊIn particular, these are minimal generators if and 
only if $F$ is not a minimal generator of $I$ (Lemma \ref{BDL lemma} 
(e)). ÊIn this case the graded Betti numbers are uniquely determined.

Hence we have to see when it can happen that the $F$ chosen in any step may 
(or may not) be a minimal generator. ÊThe point is that we are constructing 
$\Z$ inductively. ÊAt each step we are adding some set of points on a line, 
or some set of points on two lines. ÊThe graded Betti numbers are uniquely 
determined if, for regularity or other reasons, the number of points to be 
added forces $F$ to have a degree such that $F$ has no chance to be a 
minimal generator of $I$ (e.g.\ the degree is too large). ÊAlternatively if 
there is no such prohibition, we have to show that some choices of the 
points to be added correspond to $F$ a minimal generator of $I$, and other 
choices of the points to be added correspond to $F$ not a minimal generator 
of $I$.

By mimicking Example \ref{not unique} (see also Proposition \ref{std pseudo 
has std os}) we see that the standard pseudo linear configuration always 
gives an example where $F$ is not a minimal generator of $I$. ÊTherefore, 
to prove (a) we have to show that the given condition is equivalent to the 
condition that at each step, $F$ is {\em forced} to {\em not} be a minimal 
generator. ÊNotice that if at any step there is a choice between choosing 
$F$ a minimal generator or not, then not only are the graded Betti numbers 
at that step not uniquely determined, but neither are the graded Betti 
numbers for any subsequent step.

Assume first that $\Delta T'$ contains no subsequence in the list (\ref{bad 
list}). Ê Abusing notation slightly, suppose that at some intermediate step 
we have a pseudo linear configuration $\Z$ that has been obtained from the 
previous step $\X$ by a basic double link, using $F \in I_\X$ and thus 
adding a set $\Y$ to $\X$ to obtain $\Z$. ÊIf this basic double link uses a 
linear form then it corresponds to a single entry in $\Delta T'$; if it 
uses a quadric then it corresponds to a subsequence $(b,0)$ in $\Delta 
T'$. ÊWe will assume inductively that the graded Betti numbers of $I_\X$ 
are uniquely determined, and
see that then the hypothesis forces that of $\Z$ to also be uniquely 
determined.

If this basic double link corresponds to a single entry in $\Delta T'$ 
which is 2 or greater, then by Theorem \ref{pseudo k con = bdl}, $\deg F$ 
is greater than the regularity of $I_\X$ so $F$ cannot be a minimal 
generator of $I_\X$. ÊSuppose that this basic double link corresponds to a 
single entry in $\Delta T'$ which is 1. ÊIf what precedes this 1 is not a 
sequence $(\dots, 0,1,1,\dots, 1)$ then again $\deg F$ is greater than the 
regularity of $I_\X$ by Theorem \ref{pseudo k con = bdl}, so $F$ cannot be 
a minimal generator.

Next, suppose that this basic double link corresponds to an entry in 
$\Delta T'$ which is a 1, and that in $\Delta T'$ it is preceeded by 
$(\dots,b, 0, 1, 1, \dots, 1)$ (where the number of 1's may be zero). ÊBy 
the definition of a pseudo linear configuration, $b \neq 0$. ÊBy 
hypothesis, $b \neq 1$, and if $b=2$ then it is not preceeded by any 
sequence $(1,0)$, $(1,0,2,0)$, etc. ÊWe will analyze the cases $b \geq 3$ 
and $b=2$ separately, but first we make some general observations.

Corresponding to the subsequence of $\Delta T'$ given by $(\dots, b, 0, 1, 
1, \dots, 1)$, consider the sequence of configurations
\[
\dots, \X_1, \X_2,\X_3,\dots,\X_\ell = \X, \Z.
\]
Here $\X_1$ is the configuration obtained prior to this subsequence, i.e.\ it
corresponds to the initial dots before $b$ in $\Delta T'$. ÊSuppose that 
the maximum number of collinear points on $\X_1$ is $m$. Ê$\X_2$ is then 
obtained from $\X_1$ by adding two sets of collinear points, each of which 
contains $m+b$ ($\geq m+2$) points. Ê(This corresponds to the $(b,0)$ in 
$\Delta T'$.) Ê Translating to basic double links, $\X_2$ is obtained from 
$\X_1$ by a basic double link using a quadric, $Q$, and a form $F_1 \in 
I_{\X_1}$. ÊEach
subsequent Êbasic double link uses a linear form.

Note that $F_1$ is not a minimal generator of $I_{\X_1}$ (because of the 
regularity and $b \geq 2$). Ê Suppose that the minimal generators of 
$I_{\X_1}$ are $G_1,\dots,G_r$ and the graded Betti numbers of $\X_1$ are 
uniquely determined by the type (by induction). ÊThen
\[
\begin{array}{rcllccccccc}
I_{\X_2} & = & (Q G_1,\dots,Q G_r , F_1) & \hbox{where $F_1 \in I_{\X_1}$, 
not a
minimal generator of $I_{\X_1}$}
\\ I_{\X_3} & = & (L Q G_1,\dots, L Q G_r, L F_1, F_2) & \hbox{where $F_2 
\in I_{\X_2},
\deg F_2 = \deg F_1 +1$}.
\end{array}
\]
Now, if $b \geq 3$ then $\deg F_1 \geq \reg I_{\X_1} +2$. ÊHence $\deg F_1 
\geq \deg QG_i$ for all $i$, and so $F_2$ (having degree $\deg F_1 +1$) 
cannot be a minimal generator of $I_{\X_2}$ and so the listed generators of 
$I_{\X_3}$ are minimal. ÊThe same trickles down to the step from $\X$ to 
$\Z$, proving that the graded Betti numbers of $\Z$ are uniquely determined.

Now suppose that $b = 2$, but it is not preceeded by any Êsequence $(1,0)$,
$(1,0,2,0)$, etc.. ÊAgain suppose that the subsequence $(b,0)=(2,0)$ 
corresponds to a basic double link $I_{\X_2} = Q I_{\X_1} + (F_1)$ as 
above, where $I_{\X_1} = (G_1,\dots,G_r)$. ÊNow the pseudo type vector 
itself has the form $(\dots,p,q, m, m+2, m+2, m+3, m+4,\dots)$, where 
$\X_1$ is a pseudo
linear configuration of pseudo type $(\dots,p,q,m)$. ÊIf $q \leq m-2$ then 
it follows immediately that $\reg I_{\X_1} = m$ and each subsequent step 
uses an
$F$ that is not a minimal generator (not necessarily from a regularity 
argument, but rather from an analysis of the ideal as above, using Lemma 
\ref{gen lemma}). ÊIf $q = m-1$ or $q=m$, then the only danger is that
$\reg I_{\X_1} = m+1$ and that furthermore $I_{\X_1}$ has a minimal 
generator $G$ of degree $m+1$, Êso that $QG \in I_{\X_2}$ is a minimal 
generator of degree $m+3$ and can be used to construct $\X_3$ (thanks to 
the above analysis). ÊThe condition that $\reg I_{\X_1} = m+1$ holds if and 
only if the first difference of the pseudo type vector for $\X_1$ ends 
either with a 0 or with a sequence $(0,1,\dots,1)$, by Theorem \ref{pseudo 
k con = bdl}.

So we are reduced to the two cases
\[
\Delta T' = (\dots, 0,2,0,1,1,\dots,1) \hbox{ or } \Delta
T' = (\dots, 0,1,\dots,1,2,0,1,1,\dots, 1).
\]
In these cases, when can it happen that $\X_1$ has a minimal generator of 
degree
$m+1$? ÊA little thought using Lemma \ref{gen lemma} shows that in either 
case it
requires that the first 0 be preceeded by a 1, a $(1,0,2)$, a 
$(1,0,2,0,2)$, etc. ÊBut these are eliminated by our hypotheses.

Conversely, suppose that $\Delta T'$ does contain one of the subsequences 
$(1,0,1)$, $(1,0,2,0,1)$, $(1,0,2,0,2,0,1)$, etc. ÊWe know that it is 
possible to carry out the basic double links using polynomials $F$ at each 
step that are not minimal generators (mimicking Example \ref{not 
unique}). ÊSo to show non-uniqueness of the graded Betti numbers we have to 
show that at least once it is possible to choose $F$ to be a minimal 
generator in these cases.

First we consider the case where $\Delta T'$ contains a subsequence 
$(1,0,1)$. ÊHence $T'$ contains a subsequence $(m, m+1,m+1,m+2)$. ÊConsider 
a sequence of pseudo linear configurations $\X_1,\X_2,\X_3$ where
\[
\begin{array}{rcll}
I_{\X_2} & = & Q I_{\X_1} + (F_1) Ê& \hbox{where $F_1 \in I_{\X_1}, \deg 
F_1 = m+1$,}
\\ I_{\X_3} & = & L I_{\X_2} + (F_2) & \hbox{where $F_2 \in I_{\X_2}, \deg 
F_2 = m+2$}.
\end{array}
\]
The construction of basic double linkage guarantees that $I_{\X_1}$ has a 
minimal
generator of degree $m$. Ê(Notice that it cannot have a minimal generator 
of degree $m+1$ because if it did, the regularity of $I_{\X_1}$ would be 
$m+1$, so $\Delta T'$ would have a subsequence $(0,1,\dots,1,0)$, violating 
(\ref{condition}). ÊHence $F_1$ cannot be a minimal generator of 
$I_{\X_1}$.) But then
$I_{\X_2}$ has a minimal generator of degree $m+2$. ÊHence $F_2$ can either 
be chosen to be a minimal generator, or not (as illustrated in Example 
\ref{not unique}).

The analysis for the case when $\Delta T'$ has one of the other subsequences
$(1,0,2,0,1)$, $(1,0,2,0,2,0,1)$, etc.\ is very similar and is left to the 
reader.

For (b) and (c), the condition that $\Delta T'$ contains none of these 
subsequences means (according to the proof of (a)) that each basic double 
link adds a new generator. ÊAn entry of 0 in $\Delta T'$ corresponds to a 
repetition in $T'$, which in turn corresponds to the fact that two entries 
of $T'$ come from a single basic double link. ÊThe result follows immediately.
\end{proof}

A pseudo linear configuration of Êtype $T' = (\dots, m, m \dots)$ 
satisfying (\ref{condition}) can be viewed as being obtained from a pseudo 
linear configuration of type $(\dots, m, m+1)$ by removing a 
point. ÊIndeed, the only danger is that $T'$ included $(\dots,m, 
m,m+1,m+1,\dots)$ (since putting the point back would give a
configuration that has three lines with $m+1$ points), but this violates
(\ref{condition}).

In particular, let $\X$ be a linear configuration in $\P^{2}$ of type $T$, and
remove one point $P\in \X$, giving Êa pseudo linear configuration (which is
possibly still in fact a linear configuration) of type $T'$ obtained in the 
obvious way. ÊFrom our results above we recover the fact (cf.\ 
\cite{sindipaper}) that the Hilbert function of $\X \setminus P$ is 
determined from the line that $P$ lies on, but in fact we get more:

\begin{corollary} \label{remove pt}
If $\X$ is a linear configuration in $\P^2$ and $P \in \X$ then the graded 
Betti
numbers of $\X \setminus P$ are determined from the line that $P$ lies on.
\end{corollary}

\begin{proof}
When we remove a point of $\X$, we obtain a type $T'$ that has at most one 
repetition, i.e.\ $\Delta T'$ contains at most one zero. ÊHence the only 
danger, according to Theorem \ref{betti unique}, is that $\Delta T'$ 
contain the subsequence $(1,0,1)$. ÊHowever, this means that $T'$ contains 
the subsequence $(m, m+1, m+1, m+2)$ for some $m$, and it is clear that 
adding a point to any line will not give an allowable type vector for a 
linear configuration.
\end{proof}

\begin{remark}
Note that this result is not true for arbitrary $k$-configurations since 
even a $k$-configuration of type $(2,3)$ that is not a linear 
configurations provides a counterexample. ÊIt would be interesting to 
determine to what extent the result of Corollary\ref{remove pt} extends to 
$\P^n$.
\end{remark}


\section{First applications to double point schemes}

As remarked earlier, linear configurations have the property that their 
type completely determines their Hilbert function and graded Betti numbers, 
and these latter are maximal among all zero-dimensional schemes with the 
same Hilbert function. ÊIn this section and the next we are interested in 
seeing to what extent these properties are preserved for sets of 2-fat 
points which are Êsupported on a linear configuration, i.e.\ for the first 
infinitesimal neighborhood of a linear configuration.

We will use the machinery of pseudo linear configurations as an important 
component of our study, and indeed the heart of this material is the 
observation that there are surprisingly few differences between these two 
situations! In this section our focus will be to find the analog of standard
pseudo linear configurations for double point schemes supported on linear 
configurations. ÊThe key point will be that such a double point scheme, 
$\Z$, can be constructed (by basic double linkage) starting with an 
arbitrary type, $T$, by choosing the underlying linear configuration $\X$ 
in a suitable way, much as the standard pseudo linear configuration was 
chosen (for an arbitrary pseudo type vector) in a suitable way. ÊThe idea 
will be to pass to the
pseudo type vector associated to $T$ (see below). ÊThis will lead to the 
conclusion, analogous to Proposition \ref{std pseudo has std os}, that for 
any type $T$ there is a linear configuration of type $T$ whose 
corresponding double point scheme has Hilbert function computed by the 
O-sequence computation in Definition \ref{standard os}.

\begin{definition}
Let $\bar \X$ be a linear configuration of type $T = (n_1,\dots,n_r)$. ÊThe 
{\em
associated pseudo type vector} of $\bar \X$ Ê(or of $T$) is the vector $T' 
= (n_1, 2n_1, n_2, 2n_2, \dots, n_r, 2n_r)^{ord}$, where $(\ \ \ )^{ord}$ 
means that we list the entries in non-decreasing order. ÊNote that $T'$ is 
in fact a pseudo type vector, since $n_i < n_{i+1}$ for all $i$, so at most 
two entries of $T'$ take any particular value (and that happens if and only 
if we have $n_i = 2n_j$ for some $i$ and $j$).
\end{definition}

\begin{example} \label{another ex}
We will be using the associated pseudo type vector of the linear 
configuration $\X$ to build up a collection of 2-fat points with support 
$\X$. ÊWe illustrate the way we will do this with an example.

Let $\X$ be a linear configuration of type $(2,3,4)$. ÊThis gives a pseudo 
type vector of type $(2,3,4,4,6,8)$. Ê ÊWe will build up the 2-fat point 
scheme with support $\X$ (using a sequence of basic double links) in 5 steps.

\begin{itemize}
\item[Step 1:] Choose the 2 points of $\X$ on the first line. Ê(This uses 
the ``2'' in the pseudo type vector.)

\item[Step 2:] Form a basic double link to produce the scheme which 
consists of the 2 points on the first line of $\X$ and the 3 points on the 
second line. Ê(This uses the ``3'' in the pseudo type vector.)

\item[Step 3:] Form a basic double link on the ideal of Step 2 to 
(simultaneously) fatten up the two points on line 1 and add the four points 
on line 3 to the previous scheme. Ê(This uses the ``4,4'' in the pseudo 
type vector.)

\item[Step 4:] Form a basic double link on the ideal of Step 3 to fatten up 
the three points on the second line of $\X$. Ê(This uses the ``6'' in the 
pseudo type vector.)

\item[Step 5:] Form a basic double link on the ideal of Step 4 to fatten up 
the four points on the third line of $\X$. Ê(This uses the ``8'' of the 
pseudo type vector.)
\end{itemize}
The justification for why these steps are possible will be different in 
this section and the next. ÊIn this section, much as in Proposition 
\ref{std pseudo has std os}, it will be clear because of the geometry of 
the configuration. ÊIn the next section, as in the preceeding section, it 
will come as a result of showing that numerical conditions force 
conclusions about the regularity that guarantee the result. \qed
\end{example}

In the next section we will discuss the 2-type vectors, $T$, that have the 
property that {\em every} linear configuration, $\X$, of type $T$ has the 
property that its first infinitesimal neighborhood has the Hilbert function 
and graded Betti numbers described in Corollary~\ref{fat = bdl}; that is, 
for such 2-type vectors, the Hilbert function and Êgraded Betti numbers of 
any set of double points with such a support are uniquely determined.

In this section, though, we give a construction that gives, for {\em any} 
2-type vector $T$, an explicit saturated ideal of double points whose 
support is a (particular) linear configuration of type $T$. ÊWe will also 
see that sometimes there can be more than one Hilbert function for double 
points whose support is a linear configuration of type $T$. In the next 
section we will describe exactly when this happens. To illustrate the 
ideas, we begin with an example.

\begin{example} \label{ex 2 4 5}
Let $T = (2,4,5)$ be a 2-type vector, so the associated pseudo type vector is
$T' = (2,4,4,5,8, 10)$. ÊNote that we will eventually show 
(Theorem~\ref{fat = bdl}) that a double point scheme supported on {\bf any} 
linear configuration of type $T$ has the same Hilbert function but not 
necessarily the same graded Betti numbers. ÊThis is not why we have chosen 
this example. ÊThis example was selected because of its simplicity and the 
fact that it illustrates the method used.

We will form a special linear configuration, namely the ``spread out'' 
configuration (placing the points on suitable integer lattice points in the 
plane -- see Definition \ref{spread out config}), but we will also place 
``imaginary'' points to properly position the
points in which we are interested. ÊIn this case, we get the following:
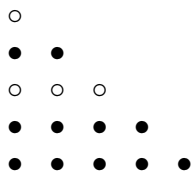
\[
\begin{array}{cccccccccccc}
\circ \\
\bullet & \bullet \\
\circ & \circ & \circ \\
\bullet & \bullet & \bullet & \bullet \\
\bullet & \bullet & \bullet & \bullet & \bullet
\end{array}
\]
We also consider three families of ``parallel'' lines: $\{L_1,L_2, 
L_3,\dots \}$,
$\{M_1,M_2,M_3,\dots \}$, and $\{D_1, D_2,D_3,\dots \}$, as follows.

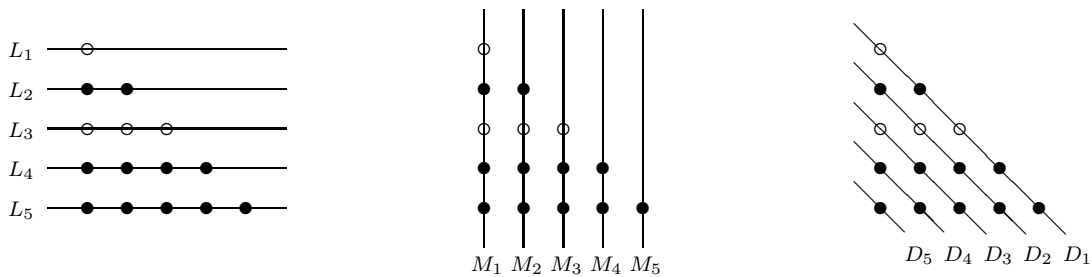
\newsavebox{\build}
\savebox{\build}(200,75)[tl]
{
\begin{picture}(200,125)
\put (-15,8){\line (1,0){90}}
\put (-15,23){\line (1,0){90}}
\put (-15,38){\line (1,0){90}}
\put (-15,53){\line (1,0){90}}
\put (-15,68){\line (1,0){90}}
\put (-30,5){$\scriptstyle L_5$}
\put (-30,20){$\scriptstyle L_4$}
\put (-30,35){$\scriptstyle L_3$}
\put (-30,50){$\scriptstyle L_2$}
\put (-30,65){$\scriptstyle L_1$}
\put (-3,5){$\bullet$}
\put (12,5){$\bullet$}
\put (27,5){$\bullet$}
\put (42,5){$\bullet$}
\put (57,5){$\bullet$}
\put (-3,20){$\bullet$}
\put (12,20){$\bullet$}
\put (27,20){$\bullet$}
\put (42,20){$\bullet$}
\put (-3,35){$\circ$}
\put (12,35){$\circ$}
\put (27,35){$\circ$}
\put (-3,50){$\bullet$}
\put (12,50){$\bullet$}
\put (-3,65){$\circ$}
\end{picture}
}

\newsavebox{\buildd}
\savebox{\buildd}(200,75)[tl]
{
\begin{picture}(200,125)
\put (0,-7){\line (0,1){90}}
\put (15,-7){\line (0,1){90}}
\put (30,-7){\line (0,1){90}}
\put (45,-7){\line (0,1){90}}
\put (60,-7){\line (0,1){90}}
\put (-5,-16){$\scriptstyle M_1$}
\put (10,-16){$\scriptstyle M_2$}
\put (25,-16){$\scriptstyle M_3$}
\put (40,-16){$\scriptstyle M_4$}
\put (55,-16){$\scriptstyle M_5$}
\put (-3,5){$\bullet$}
\put (12,5){$\bullet$}
\put (27,5){$\bullet$}
\put (42,5){$\bullet$}
\put (57,5){$\bullet$}
\put (-3,20){$\bullet$}
\put (12,20){$\bullet$}
\put (27,20){$\bullet$}
\put (42,20){$\bullet$}
\put (-3,35){$\circ$}
\put (12,35){$\circ$}
\put (27,35){$\circ$}
\put (-3,50){$\bullet$}
\put (12,50){$\bullet$}
\put (-3,65){$\circ$}
\end{picture}
}

\newsavebox{\builddd}
\savebox{\builddd}(200,75)[tl]
{
\begin{picture}(200,125)
\put (-10,78){\line (1,-1){80}}
\put (-10,63){\line (1, -1){65}}
\put (-10,48){\line (1,-1){50}}
\put (-10,33){\line (1,-1){34}}
\put (-10,18){\line (1,-1){19}}
\put (9,-12){$\scriptstyle D_5$}
\put (24,-12){$\scriptstyle D_4$}
\put (39,-12){$\scriptstyle D_3$}
\put (54,-12){$\scriptstyle D_2$}
\put (69,-12){$\scriptstyle D_1$}
\put (-3,5){$\bullet$}
\put (12,5){$\bullet$}
\put (27,5){$\bullet$}
\put (42,5){$\bullet$}
\put (57,5){$\bullet$}
\put (-3,20){$\bullet$}
\put (12,20){$\bullet$}
\put (27,20){$\bullet$}
\put (42,20){$\bullet$}
\put (-3,35){$\circ$}
\put (12,35){$\circ$}
\put (27,35){$\circ$}
\put (-3,50){$\bullet$}
\put (12,50){$\bullet$}
\put (-3,65){$\circ$}
\end{picture}
}

\begin{picture}(160,65)
\put(30,0){\usebox{\build}}
\put(180,0){\usebox{\buildd}}
\put(330,0){\usebox{\builddd}}
\end{picture}

\vskip 1.2in

Our basic double links will be of the form $L_i \cdot I + (F)$, where $F$ is a
suitable product of the $M_i$ and $D_i$. ÊAs in the previous section, we 
will ``add rows'' (which sometimes means fattening up simple points) 
according to the dictates of the pseudo type vector. ÊIn this case, we 
begin with two simple points (at the top), which we consider as the 
complete intersection Ê$I_1 = (M_1D_1, L_2)$. ÊWe then simultaneously add 
two 4's: one will fatten up $I_1$, while the other adds four simple points 
on the fourth line. ÊThis is done by forming the ideal $I_2 = L_2 L_4 \cdot 
I_1 + (M_1 D_1 M_2 D_2)$. ÊNote that $M_1D_1M_2D_2$ is double at the two 
points on the second line, and simple at the four points of the fourth 
line, so basic double linkage does indeed do the required task: $I_2$ is 
the saturated ideal of the scheme that consists of two double points on 
$L_2$ and four simple points on $L_4$. ÊWe then form
\[
\begin{array}{rcl}
I_3 & = & L_5 \cdot I_2 + (M_1 D_1 M_2 D_2 M_3) \\
I_4 & = & L_4 \cdot I_3 + (M_1 D_1 M_2 D_2 M_3 D_3 M_4 D_4) \\
I_5 & = & L_5 \cdot I_4 + (M_1 D_1 M_2 D_2 M_3 D_3 M_4 D_4 M_5 D_5)
\end{array}
\]
Notice that at each stage, the polynomial playing the role of $F$ (the 
product of the $M_i$ and $D_i$) is a multiple of the previous one, so it is 
in the previous ideal. ÊAlso, it has no component in common with the 
polynomial $L_i$ or $L_iL_j$, so basic double linkage applies. ÊFinally, it 
gives simple points when that is called for, and double points when that is 
needed, as was argued already in Example \ref{build fat} (and will be 
formalized below). ÊThe end result is the desired configuration of double 
points.

The Hilbert function of the set of points that we have constructed is again 
obtained from a simple computation (as are the graded Betti numbers). ÊThe 
Hilbert function is

\bigskip

\begin{tabular}{ccccccccccccccccccccccccc}
\begin{tabular}{cccccccccccccccccccccccccccc}
ÊÊ & & & & & 1 & 1 \\
ÊÊ & & & 1 & 2 & 2 & 2 & 1 \\
ÊÊ & & 1 & 1 & 1 & 1 & 1 \\
ÊÊ & 1 & 1 & 1 & 1 & 1 & 1 & 1 & 1 \\
1 & 1 & 1 & 1 & 1 & 1 & 1 & 1 & 1 & 1 \\ \hline
1 & 2 & 3 & 4 & 5 & 6 & 6 & Ê3 & 2 & 1
\end{tabular}
&
$\leadsto$ & 1 & 3 & 6 & 10 & 15 & 21 & 27 & 30 & 32 & 33 & 33 & \dots
\end{tabular}

\bigskip

\noindent and the graded Betti numbers are given by the {\tt macaulay} diagram

\newpage

\begin{verbatim}
ÊÊ Ê Ê Ê Ê Ê Ê Ê Ê Ê Ê Êtotal: Ê Ê Ê1 Ê Ê 6 Ê Ê 5
ÊÊ Ê Ê Ê Ê Ê Ê Ê Ê Ê Ê Ê--------------------------
ÊÊ Ê Ê Ê Ê Ê Ê Ê Ê Ê Ê Ê Ê Ê0: Ê Ê Ê1 Ê Ê - Ê Ê -
ÊÊ Ê Ê Ê Ê Ê Ê Ê Ê Ê Ê Ê Ê Ê1: Ê Ê Ê- Ê Ê - Ê Ê -
ÊÊ Ê Ê Ê Ê Ê Ê Ê Ê Ê Ê Ê Ê Ê2: Ê Ê Ê- Ê Ê - Ê Ê -
ÊÊ Ê Ê Ê Ê Ê Ê Ê Ê Ê Ê Ê Ê Ê3: Ê Ê Ê- Ê Ê - Ê Ê -
ÊÊ Ê Ê Ê Ê Ê Ê Ê Ê Ê Ê Ê Ê Ê4: Ê Ê Ê- Ê Ê - Ê Ê -
ÊÊ Ê Ê Ê Ê Ê Ê Ê Ê Ê Ê Ê Ê Ê5: Ê Ê Ê- Ê Ê 1 Ê Ê -
ÊÊ Ê Ê Ê Ê Ê Ê Ê Ê Ê Ê Ê Ê Ê6: Ê Ê Ê- Ê Ê 3 Ê Ê 2
ÊÊ Ê Ê Ê Ê Ê Ê Ê Ê Ê Ê Ê Ê Ê7: Ê Ê Ê- Ê Ê - Ê Ê 1
ÊÊ Ê Ê Ê Ê Ê Ê Ê Ê Ê Ê Ê Ê Ê8: Ê Ê Ê- Ê Ê 1 Ê Ê 1
ÊÊ Ê Ê Ê Ê Ê Ê Ê Ê Ê Ê Ê Ê Ê9: Ê Ê Ê- Ê Ê 1 Ê Ê 1
\end{verbatim}

\bigskip

\noindent This fails to be the graded Betti numbers of the lexsegment ideal 
with the given Hilbert function only because one step involved a basic 
double link using a quadric instead of a linear form (to produce $I_2$) -- 
see Theorem \ref{betti unique}. \qed
\end{example}

With this example giving the reader our basic ideas, we are now ready to 
extend Proposition \ref{std pseudo has std os} to 2-fat points. ÊNotice that
part $iii)$ of the following theorem is much cleaner than parts $iii)$ and 
$iv)$ of Proposition \ref{std pseudo has std os}, because of the 
extra Ê``compactness" provided by the non-reducedness step. ÊNote also that 
simply using the standard lifting, without ``raising'' the rows to fit into 
the isosceles triangle, is not enough. ÊSee Example \ref{special}.

\begin{theorem}\label{any 2-type}
Let $T = (n_1,\dots,n_r)$ be a 2-type vector, with $n_1 < n_2 < \dots < 
n_r$, and let $T' = (m_1,\dots,m_{2r})$ be the associated pseudo type 
vector. ÊLet
$\X$ be the spread out linear configuration with type vector $T$ (see 
Definition
\ref{spread out config}). ÊLet $\Z$ be the set of 2-fat points supported on 
$\X$. ÊThen

\begin{itemize}
\item[i)] $\Z$ can be built up by basic double linkage.

\item[ii)] The first difference of the Hilbert function of $\Z$ is the standard
O-sequence associated to $T'$ (from the O-sequence computation in 
Definition \ref{standard os}).

\item[iii)] ÊThe regularity of $\Z$ is $m_{2r} = 2n_r$.

\end{itemize}
\end{theorem}

\begin{proof}
We have $T' = (n_1, 2n_1, n_2, 2n_2, \dots, n_r, 2n_r)^{ord}$. ÊNote that 
if an entry of $T'$ is odd then it only occurs once, and if an entry is 
even then it occurs {\em at most} twice. ÊWe let $m_i$ denote the (ordered) 
entries of $T'$, so $m_1 = n_1, \dots, m_{2r} = 2n_r$.

Consider the spread out linear configuration, $\X$, with Ê$r$ rows, each 
having $n_r$ points, as in Example~\ref{ex 2 4 5}. Ê Again as in 
Example~\ref{ex 2 4 5}, we consider three families of lines, 
$\{L_1,\dots,L_r\}$, $\{M_1, \dots, M_{n_r}\}$, $\{D_1,\dots,D_{n_r} 
\}$. ÊThe $L_i$ are the ``horizontal'' lines, the $M_i$ are the 
``vertical'' lines and the $D_i$ are the ``diagonal'' lines (starting at 
the ``hypotenuse'').

A basic double link has the form $I_{t+1} = G \cdot I_t + (F)$, where $F 
\in I_t$ and $(F,G)$ is a regular sequence. In our case, at each step the 
role of $F$ will be played by a suitable product $M_1D_1M_2 D_2 \cdots$, 
alternating between them. ÊIf we have completed the construction for 
$m_{i-1}$ in the pseudo type vector $T'$, then to build the ideal 
corresponding to the entry $m_i$ (which may or may not be equal to 
$m_{i+1}$) in the pseudo type vector, Êthe number of factors in the 
polynomial $F$ is equal to $m_i$. ÊEach subsequent $F$ will build on the 
ones before by adding factors consisting of the $M_i$ and $D_i$. ÊThis 
guarantees that we will always have $F \in I_t$, and in fact by checking 
the regularity we can see that the $F$ we are using is not a minimal 
generator of $I_t$.

The role of $G$ will always be played by either one $L_j$ (if $m_i < 
m_{i+1}$), or a product of two $L_j$ (if $m_i = m_{i+1}$), as dictated by 
$T'$.

We make the following observations:

\begin{enumerate}
\item $m_1 = n_1 < m_2$. ÊThe construction starts with the ideal $I_1$ that 
is the complete intersection of $L_1$ and $F$, where $F$ is the product 
$M_1 D_1,\dots$, taking $n_1$ factors. ÊThis is the ideal of $n_1$ simple 
points on $L_1$.

\item \label{setup 1=} At any step, if $m_i < m_{i+1}$ then $F$ has $m_i$ 
factors, and either it cuts out $m_i$ simple points on $L_j$ or else we 
have $m_i = 2n_k$ for some $k < i$, and $F$ is double at each of $n_k$ 
points on $L_k$.

\item If $m_i < m_{i+1}$, $I_t$ is the current ideal, and $F$ is chosen as in
(\ref{setup 1=}), then the ideal $I_{t+1} = L_j \cdot I_t + (F)$ is a 
saturated ideal that either adds $n_j$ simple points on $L_j$ or else it 
``fattens up" (doubles) $n_k$ points on $L_k$, respectively.

\item \label{setup 2 =} At any step, if $m_i = m_{i+1}$ then one of them 
(without loss of generality say it is $m_i$) is the term $2n_j$ for some $j 
< i$, and the other is equal to $n_k$ for some $k \leq i$. ÊIn this case 
$F$ has $m_i$ factors, and it has $n_j$ singular (double) points along the 
line $L_j$ ($1 \leq j \leq r$) and $n_k$ simple points along the line $L_k$ 
($1 \leq k \leq r$).

\item If $m_i = m_{i+1}$, $I_t$ is the current ideal, and $F$ is chosen as in
(\ref{setup 2 =}), then the ideal $I_{t+1} = L_j L_k \cdot I_t + (F)$ is a 
saturated ideal that adds $n_k$ simple points on $L_k$ and ``fattens up'' 
$n_j$ already-existing simple points on $L_j$.
\end{enumerate}

The end result, after completing this procedure by reaching $m_{2r}$, is the
saturated ideal of double points supported on the spread out linear 
configuration of type $T$. ÊThe computation of the Hilbert function is 
identical to that in Theorem~\ref{pseudo k con = bdl}. ÊThis completes $i)$ 
and $ii)$.

Now, the numerical information obtained from the basic double linkage is 
identical to that we saw in the reduced situation -- it only depends on the 
degrees of the polynomials, and not on the geometry of the 
singularities. ÊIn particular, we obtain from Proposition \ref{std pseudo 
has std os} Êthat what can play havoc with the regularity here is the 
existence of certain subsequences in $\Delta T'$. ÊIn particular, if 
$\Delta T'$ has an entry that is $>1$ between any two zero entries then the 
regularity can only be $m_{2r}+1$ or $m_{2r}$, depending (respectively) on 
whether $\Delta T'$ ends with one of the subsequences 0, $(0,1)$, 
$(0,1,1)$, \dots, or not. ÊIf $\Delta T'$ has zero entries between which 
there are only 1's then the regularity can (in principle) be arbitrarily 
bigger than $m_{2r}$. ÊSo we have to verify that such things
cannot happen for 2-fat points.

First note that $\Delta T'$ can {\em not} end with a 0 or a 1. ÊIndeed, we have
$m_{2r} = 2n_r$, which is even, and if $m_{2r-1} = 2n_r$ or $2n_r -1$ then 
this entry is not the double of a previous one and hence its double is 
still to come. ÊSo if, between any two zero entries of $\Delta T'$, there 
is at least one entry $>1$, we now know that the regularity of $\Z$ is 
$m_{2r} = 2n_r$.

It is certainly possible for $\Delta T'$ to have Êa subsequence 
$0,1,\dots,1,0$. ÊFor instance, take $T = (8,9,10,16,17,19,20)$; then 
$\Delta T' = (8,
1,1,6,0,1,1,1,1,0,12,2,4,2).$ ÊIt is clear from the discussion leading to the
O-sequence computation in Definition \ref{standard os} that at each step, 
if we are performing (without loss of generality) a basic double link with 
$\deg F = m$ (say) and $\deg G = 1$, building from a zero-dimensional 
scheme $\X$ to $\Y$, then
\[
\reg(\Y) = \max \{ \deg F, \reg(\X) +1\}.
\]
The point that we will make is that (as we have seen) what creates 
``problems" for the regularity is a double occurrence of an integer in 
$T'$, say $(\dots,m,m,\dots)$, i.e.\ a 0 in $\Delta T'$. ÊBut such an 
occurrence automatically forces a $2m$ also in $\T'$, and this corrects the 
problems.

Indeed, suppose that the {\em last} 0 in $\Delta T'$ occurs in position $d$,
and that prior to this 0 there are $k$ zeros. ÊSo $ T' = (\dots, 
m,m,\dots)$, where the second $m$ occurs in position $d$. Clearly $k \leq 
\frac{m}{2}$ (a zero in $\Delta T'$ has to correspond to an even number in 
$T'$). Then the regularity of the subscheme produced up to that point in 
$T'$ is $\leq m+k \leq m + \frac{m}{2}$. ÊWhat can happen after this point 
in $T'$? ÊEither all of the remaining entries are of the form $2n_i$ (so 
the last one is $2m =2n_r$), or there are more $n_i > m$ (so the last entry 
of $T'$ is $2n_r > 2m$).

In the first case, the number of remaining steps is clearly $\leq 
\frac{m}{2}$, since the number of remaining steps is exactly the number of 
$n_i$ for which $2n_i > m$. ÊHence the regularity of the final double point 
scheme is
\[
\max \left \{ 2m, \left ( m+\frac{m}{2} \right ) +\frac{m}{2} Ê\right \} = 
2m = 2n_r.
\]

In the second case, when we reach the entry $2m$ in $T'$, we already have 
regularity being determined by the entry (namely $2m$ in this case), and 
each subsequent entry preserves this property. ÊHence again the regularity 
of the resulting scheme is $2n_r$.
\end{proof}

\begin{example} \label{special}
ÊÊThe construction in this section sometimes has very special 
properties. ÊFor example, suppose that we want to study the Hilbert 
function of the first infinitesimal neighborhood of a linear configuration 
of type $T = (4, 5, 8, 9, 10)$. ÊThe basic double link prediction for this 
Hilbert function is
$$
1, 2, 3, 4, 5,6 ,7, 8, 9, 10, 10, 10, 10, 7, 4, 3, 3, 3, 2, 1.
$$
But even the standard lifting of the lex-segment ideal (putting the points 
on the
integer lattice points) gives the more general Hilbert function
$$
1, 2, 3, 4, 5,6 ,7 ,8, 9, 10, 10, 10, 10, 8, 3, 3, 3, 3, 2, 1.
$$
(which is also the Hilbert function for the first infinitesimal 
neighborhood of a
generically chosen linear configuration of this type). ÊBut moving the 
points ``upward'' as indicated in this section, to add collinearity of the 
``diagonal'' points, is enough to change the value in degree 13 to this 
more special function. 
ÊWe have verified this on {\tt macaulay} \cite{macaulay}. Notice
that the Êassociated pseudo type vector $T'$ does  not satisfy (\ref{condition}). \qed
\end{example}

\begin{remark} It may be noted that a key difference between Theorem 
\ref{any 2-type} and Theorem \ref{pseudo k con = bdl} is that in the latter 
we had to use vanishing of first cohomology to guarantee lifting of 
non-zero elements, which in Theorem \ref{any 2-type} is not guaranteed 
simply by the cohomology; rather, we used the simplicity of the geometry to 
guarantee the existence of suitable curves (the unions of the lines). \qed
\end{remark}

\begin{remark}
The construction of Theorem~\ref{any 2-type} would work equally well if the 
families $\{L_1,L_2,\dots \}$, $\{M_1,M_2,\dots \}$ and $\{D_1,D_2,\dots 
\}$ (each of which has a common point at infinity) were replaced by three 
different families of lines in $\P^2$, each with a common point in $\P^2$. \qed
\end{remark}


\section{When are the Hilbert function and graded Betti numbers uniquely 
determined?}

In this section we will show how to apply the ideas of Theorem \ref{pseudo 
k con = bdl} and Theorem \ref{betti unique}, and especially their proofs, 
to the study of double points in $\P^2$. ÊWe will show that the same ideas 
in fact produce the (non-reduced) double point scheme by basic double 
linkage, and the same kind of uniqueness results continue to hold. ÊSome of 
the important ideas used here were illustrated in Example \ref{build fat}.

The following is the main result of this section, and extends to 2-fat 
points the
results on Hilbert functions and graded Betti numbers of pseudo linear 
configurations.

\begin{theorem} \label{fat = bdl}
Let $\bar \X$ be a linear configuration of type $T = (n_1,\dots,n_r)$, and 
let $T' =$ \linebreak $(m_1,\dots, m_{2r})$ be the associated pseudo type 
vector. ÊLet $\bar \Z$ be the 2-fat Êpoint scheme Êsupported on $\bar \X$.

\begin{itemize}

\item[(a)] Assume that for each $i$ we have the property (\ref{condition}) 
of Theorem \ref{pseudo k con = bdl}, namely that between any two zero 
entries of $\Delta T'$ there is at least one entry that is $>1$. ÊThen 
$\bar \Z$ can be constructed as a sequence of basic double links, and its 
Hilbert function is uniquely determined and can be computed by the 
O-sequence computation of Definition Ê\ref{standard os}.

\item[(b)] Conversely, if (\ref{condition}) does not hold then there are
linear configurations of the given type, $T$, whose corresponding double 
points do {\em not} arise by basic double linkage. ÊFurthermore there are 
two different
linear configurations of type $T$ such that the corresponding double points 
have
different Hilbert functions.

\item[(c)] Assume again that (\ref{condition}) holds. ÊAssume further that 
$\Delta T'$ contains no subsequence $(1,0,1)$, $(1,0,2,0,1)$, 
$(1,0,2,0,2,0,1)$, etc. ÊThen, in addition, the graded Betti numbers of 
$I_\Z$ are uniquely determined, as described in Theorem \ref{betti unique}.

\item[(d)] Conversely, if (\ref{condition}) holds, but $\Delta T'$ {\em 
does} contain a subsequence $(1,0,1)$, $(1,0,2,0,1)$, $(1,0,2,0,2,0,1)$, 
etc.\ then there are two different linear configurations of type $T$ such 
that the corresponding double points have the same Hilbert function (by 
part (a)), but the graded Betti numbers are different.
\end{itemize}
\end{theorem}

\begin{proof}
As usual we assume that $\bar \X=\bigcup_{i=1}^r \X_i$, where $\X_i$ 
consists of $n_i$ points on line $L_i$, for Ê$1 \leq i \leq r$. By the 
definition of
a linear configuration, $n_{i-1} < n_i$ for all $i$. ÊIf we set $L=L_i$ then
$[I_{\bar \Z} : L]$ is a saturated ideal defining the union of ${\X_i}$ 
(the reduced points on $L$) and the double points whose supports are not on 
$L$. ÊFurthermore, $\frac{I_{\bar \Z} + (L)}{(L)}$ is the (non-saturated) 
ideal of a subscheme of $L$ that has degree $2n_i$ and is supported on 
$\X_i$ with degree Êtwo at each point and tangent direction given by $L$.

Our strategy will be to consider $\bar \Z$ inductively as a ``limit'' pseudo
linear configuration of type $T'$, and to construct $\bar \Z$ in the order 
dictated by $T'$, just as in Theorem \ref{pseudo k con = bdl} (see Example 
\ref{build fat}and Example \ref{another ex}). ÊAgain, if we have reached 
and completed $m_{i-1}$ in our construction, then the next step will handle 
$m_i$ alone if $m_i < m_{i+1}$, and it will handle $m_i$ and $m_{i+1}$ 
simultaneously if $m_i = m_{i+1}$ (which then is necessarily an even 
number). ÊWhen we have $m_i = m_{i+1} = 2m_j$ for some $i$ and $j$, this 
will involve simultaneously ``fattening up'' the points corresponding to 
$m_j$ and adding the simple points corresponding to $m_i$. Ê(The next 
section applies this idea in a more Êconcrete, geometric way.) ÊNote that 
any intermediate step may or may not produce a scheme consisting entirely 
of double points. ÊOnly the final result will necessarily consist entirely 
of double points, namely $\bar \Z$. ÊNote also that $L$ does not 
necessarily progress monotonically through the $L_i$, since when it 
``fattens up'' a set of points on a line, that line will be a
previously considered one (as illustrated in Example \ref{build fat} and 
Example
\ref{another ex}).

The ``fattening up'' process is based on the following observation: if $P$ 
is a point in $\P^2$ and if $L_1, L_2, L \in I_P$ are Êlinear forms, then 
$L \cdot I_P + (L_1 L_2)$ is the saturated ideal of the double point scheme 
defined by the
(saturated) ideal $I_P^2$, as long as $L$ has no component in common with 
either $L_1$ or $L_2$. ÊMore generally, let $P$ be a reduced point of a 
scheme $\Z$, let $F \in I_\Z$ be a homogeneous polynomial such that $F \in 
I_P^2$ and $F \notin I_P^3$, and let $L \in I_P$ with no component in 
common with $F$ even locally (i.e.\ the intersection of $F$ and $L$ is a 
zero-dimensional scheme that has degree 2 at $P$). ÊThen $L \cdot I_\Z + 
(F)$ is the saturated ideal of a zero-dimensional scheme in $\P^2$, and at 
$P$ this zero-dimensional scheme is the 2-fat point supported on $P$. ÊIt 
is worth noting that if we allowed $F$ to be smooth at $P$ and $L$ were 
tangent to $F$ at $P$, then the new zero-dimensional scheme again 
would Êhave degree ($\geq$) 3 at $P$, but would be curvilinear, not ``fat."

So, mimicking the approach of Theorem \ref{pseudo k con = bdl}, suppose 
that we have reached and completed $m_{i-1}$. ÊAs before, there are two 
possibilities: either $m_i < m_{i+1}$ or $m_i = m_{i+1}$.

We first suppose that $m_i < m_{i+1}$, and we set $L$ to be the line 
containing the $m_i$ ``points'' (which is Ênot necessarily $L_i$). ÊThese 
will either be

\begin{itemize}
\item[i)] $m_i$ reduced points (which we will add singly), or

\item[ii)] $\frac{m_i}{2}$ length two schemes (not fat) on $L$, which we 
will ``add'' to $\frac{m_i}{2}$ already-existing Êsingle points to obtain
$\frac{m_i}{2}$ double points. ÊNote that then $\frac{m_i}{2}$ is one of 
the $n_j$.

\end{itemize}

\noindent In either case $\Y$ will denote this subscheme of $L$ of degree 
$m_i$, and $\X$ will denote the subscheme of $\bar \Z$ constructed 
(inductively) up to that point. Ê$\Z$ will denote the ``union" of $\X$ and 
$\Y$, but now this is more delicate to define. ÊIf $\Y$ is reduced (case 
i)), we simply take $\Z$ to be the union in the usual sense. Ê If $\Y$ is 
non-reduced, then $\Z$ will denote the scheme obtained from the scheme of 
the previous step by replacing the $\frac{m_i}{2}$ simple points with 
$\frac{m_i}{2}$ double points. ÊWe have to show that either way, $\Z$ is 
obtained from $\X$ by basic double linkage.

We again consider the exact sequence
\[
\begin{array}{cccccccccccccccccc}
0 & \rightarrow & [I_\Z : L](-1) & \stackrel{\times L}{\longrightarrow} &
I_\Z & \rightarrow & \displaystyle \frac{I_\Z + (L)}{(L)} & \rightarrow & 0. \\
&& || \\
&& I_\X (-1)
\end{array}
\]
The mechanics of the proof (using regularity to lift elements, and 
analyzing minimal generators) are identical to those of Theorem \ref{pseudo 
k con = bdl} and Theorem \ref{betti unique} and will not be repeated 
here. ÊWhat is new is the justification that it all works even in the 
non-reduced situation. ÊBut in fact, $\Y$ is a divisor on $L = \P^1$, and 
whether it is reduced or not, its ideal in $L$ begins in degree $m_i$ just 
as before. ÊA non-zero element of
$I_{\Y|L}$ (the saturation of $\frac{I_\Z + (L)}{(L)}$) in degree $m_i$ 
lifts to an element, $F$, of $(I_\Z)_{m_i}$ just as before, and $\Y$ is the 
complete intersection of $F$ and $L$. Ê We then form the ideal $I = L \cdot 
I_\X + (F)$. ÊThis is the saturated ideal of a scheme $\Z$ that is the same 
as $\X$ for points off $L$, and makes a non-reduced degree three subscheme 
of $\P^2$ at each point of the support of $\Y$. ÊThe one remaining subtlety 
is to ascertain that at each such point in the support of $\Y$, the 
non-reduced scheme that we obtain is really a 2-fat point. ÊThis would fail 
to happen, as noted above, if the polynomial $F$ is smooth at a point of 
$\Y$ and tangent to $L$ there, rather than singular there. Ê(Such an $F$ 
certainly restricts to a $\Y$
that is double at each point, as a subscheme of $L = \P^1$.) ÊBut this is
resolved by the fact that we know that we are lifting elements of $I_{\Y|L}$ to
$I_\Z$, which we knew in advance to consist of 2-fat points at each of the
$\frac{m_i}{2}$ points in the support of $\Y$. ÊHence $F$ is not smooth at 
any of those points, and must be double there. ÊSo now, $L \cdot I_\X + 
(F)$ and $I_\Z$ are both saturated ideals defining the same 
zero-dimensional subscheme, hence they are equal. ÊThis completes the proof 
of (a).

Parts (c) and (d) continue to have (\ref{condition}) as a hypothesis, 
meaning that the configurations of 2-fat points considered there 
necessarily arise by basic double linkage, but the graded Betti numbers are 
in question. ÊWe consider these parts first, and then turn to (b).

The Hilbert function and regularity of the new scheme are obtained just as 
in Theorem \ref{pseudo k con = bdl}. ÊIn case (c), the graded Betti numbers 
are produced just as in Theorem Ê\ref{betti unique}. ÊIf $m_i = m_{i+1}$, 
instead of $L$ we again use $Q$ which is the product of two linear forms. 
One of them will contain $m_i$ reduced points and the other will be viewed 
as containing $\frac{m_i}{2}$ double points, as noted above. ÊNote that 
these are distinct lines! ÊAgain the same proof as in Theorem \ref{pseudo k 
con = bdl} and
Theorem \ref{betti unique} works, with the same modifications as in the 
previous
paragraph.

For part (c), the point is that we have just shown that the double points are
constructed with liaison addition in a manner perfectly analogous to that 
used for the pseudo linear configurations. ÊThe conditions in (c) then 
guarantee that every step forces us to choose $F$ {\em not} a minimal 
generator of the previous ideal, hence the conclusion that the graded Betti 
numbers are uniquely determined.

Part (d) is slightly more subtle, however. ÊEach step of the basic double 
linkage
either adds a new set of reduced points, ``fattens up" an existing set, 
or Êdoes both simultaneously. ÊNote that there is less freedom if we are 
constrained to a previously existing support. ÊHowever, the ``fattening 
up'' process can only be done if the corresponding entry in $T'$ is 
even! ÊA subsequence $(1,0,1)$ in $\Delta T'$ corresponds to a subsequence
\[
m, m+1, m+1, m+2
\]
in $T'$, and a subsequence $(1,0,2, 0,2,\dots,0,2,0,1)$ in $\Delta T'$ with 
$k$ 2's corresponds to a subsequence
\[
m, m+1, m+1, m+3, m+3, \dots, m+2k+1, m+2k+1, m+2k+2
\]
in $T'$. ÊIn each case, the last entry must be odd ($m+2$ and $m+2k+2$, 
respectively). ÊComparing with the proof of Theorem \ref{betti unique}, it 
is exactly at this point that there is a choice of choosing $F$ a minimal 
generator or not, and since the number is odd, this must correspond to 
adding new reduced points, not ``fattening up'' already existing 
points. ÊHence we have complete freedom with $F$, and Ê(d) follows.

We now turn to (b). ÊThe proof is very similar to the last part of Theorem 
\ref{pseudo k con = bdl}, with some fine tuning. ÊWe have to show that if 
$\Delta T'$ contains a subsequence $(\dots, 0, 1, \dots, 1,0)$ (all 1's 
between the two 0's) then there exist (at least) two Êlinear configurations 
of type $T$ whose
associated double points (first infinitesimal neighborhoods) have different 
Hilbert functions. ÊWe have already noted that an associated pseudo type 
vector cannot end with a 0, so the argument will be slightly different from 
that of Theorem \ref{pseudo k con = bdl}.

Note that we showed in Theorem \ref{any 2-type} that for {\em any} Êtype
vector $T$, there always exists one Êlinear configuration (the spread out 
one) whose first infinitesimal neighborhood can be constructed by basic 
double linkage, and hence has ÊHilbert function whose first difference is 
given by the O-sequence computation of Definition \ref{standard 
os}. ÊSo Êwe have to show
that such a subsequence allows for a 2-fat point scheme that can {\em not} be
constructed entirely by basic double links, and that correspondingly the 
Hilbert
functions are different.

Suppose that we are given the type vector $T = (n_1,\dots,n_r)$, and a
linear configuration $\X$ of type $T$. ÊFrom $T$ we derive the associated 
pseudo type vector $T' = (n_1, 2n_2, \dots,n_r, 2n_r)^{ord}$. ÊThis 
information gives the recipe to ``fatten up" $\X$ to a 2-fat point scheme 
$\Z$ by basic double linkage, if such a process is possible. ÊIt is 
important to note that each entry of $T'$ corresponding to an $n_i$ 
produces $n_i$ reduced points on a new line, and each entry of $T'$ 
corresponding to a $2n_i$ ``fattens up" $n_i$ previously existing points on 
a line. ÊThe only ambiguity comes when we have two consecutive entries that 
are equal. ÊUsually we do these simultaneously, by taking $G$ to be the 
product of the two linear forms. Ê However, for this proof we will consider 
such a situation as arising from two consecutive basic double links using 
the same polynomial $F$ and taking $G$ linear, rather than one basic double 
link using $G$ quadratic.

\begin{quote} \label{convention}
{\em We will make the convention that the first basic double link 
corresponds to ``fattening up'' previously existing points, while the 
second one corresponds to producing new reduced points.}
\end{quote}

If basic double linkage is possible at each step, the end result of this 
process is the desired 2-fat point scheme $\Z$ supported on $\X$. ÊHowever, 
each intermediate step is the saturated ideal of a zero-dimensional scheme 
that is ``2-fat" at some points and reduced at others.

Suppose that $\Delta T'$ contains a subsequence $(\dots, 0,1,\dots,1,0)$ 
(all 1's
between the 0's) and consider the first occurrence of this 
subsequence. ÊThis means that
\[
T' = (m_1,\dots,m_{p-2}, m_{p-1},m_p, \dots)
\]
with Ê$2n_i = m_p = m_{p-1} = m_{p-2}+1$. Ê Let $T''$ be the pseudo type 
vector $(m_1,\dots, m_{p-2}, m_{p-1})$. ÊThen $\Delta T'"$ ends in a 
sequence $(\dots, 0, 1, \dots, 1)$ where the end consists of nothing but 
1's. Ê We have assumed that $T''$ satisfies (\ref{condition}). ÊLet $\Z_1$ 
be the zero-dimensional scheme corresponding to $T''$, following our 
procedure of basic double linkage; $\Z_1$ is supported on some subset of 
$\X$. ÊIts Hilbert function is as
described in the O-sequence computation of Definition \ref{standard 
os}. ÊFurthermore, it follows from what we have already proven that the 
regularity of $I_{\Z_1}$ is $m_{p-1}+1 = Êm_{p}+1 = 2n_i +1$. ÊThe last 
basic double link in this sequence ``fattened up" a previously existing 
$n_i$ points.

Now consider an additional line $L$, and choose $\Y$ to be a general set of 
$2n_i$ points on $L$. ÊLet $\Z_2 = \Z_1 \cup \Y$. Ê $\Z_2$ is a basic 
double link of $\Z_1$ if and only if there is a form $F \in 
(I_{\Z_1})_{2n_i}$ that
contains $\Y$ but does not vanish on $L$.

We have an exact sequence
\[
0 \rightarrow {\mathcal I}_{\Z_1} (2n_i-1) \stackrel{\times L}{\longrightarrow}
{\mathcal I}_{\Z_1}(2n_i)
\rightarrow {\mathcal O}_{L} (2n_i) \rightarrow 0.
\]
Since the regularity of $I_{\Z_1}$ is $2n_i+1$, we have
\begin{equation} \label{eexx2}
0 \rightarrow (I_{\Z_1})_{2n_i-1} \rightarrow (I_{\Z_1})_{2n_i}
\stackrel{r}{\longrightarrow} H^0({\mathcal O}_{L} (2n_i)) \rightarrow 
H^1({\mathcal
I}_{\Z_1}(2n_i-1)) \rightarrow 0
\end{equation}
where the last cohomology group is not zero. ÊChoosing $\Y$ as above is 
equivalent to choosing a general element of the vector space $H^0({\mathcal 
O}_{L}(2n_i))$. ÊThe image of $r$ is a proper subspace of $H^0({\mathcal 
O}_{L}(2n_i))$, so the general section of $H^0({\mathcal O}_{L}(2n_i))$ 
defining $\Y$ is not in the image of $r$. ÊWe conclude that any form in 
$(I_{\Z_1})_{2n_i}$ that vanishes on $\Y$ must in fact vanish on all of 
$L$. ÊHence we cannot express $\Z_2$ as a basic double link of $\Z_1$.

We claim that the value of the Hilbert function of $\Z_2$ in degree $2n_i$
differs (in fact is larger) from the value of the corresponding Hilbert 
function
given by the O-sequence computation in Definition \ref{standard os} in degree
$2n_i$ (see Remark \ref{std hf from os}). ÊIndeed, suppose that $\Z'$ were a
zero-dimensional scheme that was produced by a sequence of basic double 
linkages of the same type, and hence has the standard Hilbert function for 
that type. ÊIn degree $2n_i$, the forms that vanish on $\Z_2$ consist 
entirely of products of $L$ with forms of degree $m_p-1$ vanishing on 
$\Z_1$ (as discussed above), while $\Z'$ has those but also has a form of 
degree $2n_i$ that is not of that form. ÊHence the first ÊHilbert function 
is larger than the second in degree $2n_i$.

Now we continue along $T'$. ÊWe have reached the entry $m_p = 2n_i$ and 
constructed a zero-dimensional scheme $\Z_2$ whose Hilbert function is {\em 
not} the one predicted by the O-sequence computation of Definition 
\ref{standard os}, precisely because at the last step we added a set of 
points and showed that it could not arise by basic double linkage. ÊAt each 
subsequent step, one of three things can happen: (i) because of regularity 
arguments like those above, we are guaranteed that that step can be 
accomplished by basic double linkage; (ii) a step corresponds to 
``fattening up'' an existing set of reduced points, and it happens that it 
can be accomplished by basic double linkage, or
(iii) whether because of the position of the existing points to be 
``fattened up" or because of the free choice of general reduced points, 
basic double linkage cannot be performed.

As in Theorem \ref{pseudo k con = bdl}, if (i) or (ii) hold then the 
resulting scheme again fails to have the standard O-sequence predicted by 
the O-sequence computation of Definition \ref{standard os} because we are 
adding the expected amount to an already larger Hilbert function. ÊIn the 
third case, as in the argument just made, the Hilbert function becomes 
correspondingly larger than it would have been had basic double linkage 
been possible, hence gets even farther from the predicted O-sequence.

In the end we obtain a set of 2-fat points supported on a linear 
configuration whose Hilbert function is different from that of a set of 
2-fat points supported on a spread out configuration. ÊThis proves (b).
\end{proof}

\begin{remark}\label{same reg}
Although, as Theorem \ref{fat = bdl} states, it need not be true that the 
first infinitesimal neighborhood of two linear configurations of type 
$(n_1,\dots,n_r)$ have the same Hilbert function or the same Betti numbers 
in their minimal free resolution, there is one thing they will have in 
common -- namely their regularity (which is $2n_r$).

To see why that is so, just observe that by Lemma \ref{bd reg} the first 
infinitesimal neighborhood of any linear configuration of type $(n_1, 
\ldots , n_r)$ has regularity $\leq 2n_r$. ÊHowever, the first 
infinitesimal neighborhood also always has a subscheme of length $2n_r$ on 
a line and so the regularity is $\geq 2n_r$.
\end{remark}

\begin{example}\label{not unique betti diag}
ÊFrom Theorem \ref{fat = bdl} we see that linear configurations of the same 
type may have, for their first infinitesimal neighborhoods, the same 
Hilbert function but not the same graded Betti numbers. Ê It would be 
interesting to know exactly what the possibilities are for the Betti 
numbers of these double point schemes in such a case. ÊThis example deals 
with that situation.

Let $T = (2,3,4,5)$. ÊLet $\bar \X$ be a linear configuration of type $T$ 
and let $\bar \Z$ be the first infinitesimal neighborhood of $\bar 
\X$. ÊThen the Hilbert function of $\bar \Z$ is uniquely determined and has 
first difference
\[
1,\ 2,\ 3,\ 4,\ 5,\ 6,\ 7,\ 8,\ 5,\ 1.
\]
However, the Êgraded Betti numbers are not uniquely determined. ÊThe associated
pseudo type vector is $(2,3,4,4,5,6,8,10)$. ÊOne can check that there are 
actually {\em two} times, in making the construction of Theorem \ref{fat = 
bdl}, when there is apparently a choice between using a minimal generator 
or not, namely when we deal with the 5 and when we deal with the 
6. ÊHowever, notice that while there is freedom in choosing where the 5 
points are located, there is no such freedom for the 6 since it represents 
the ``fattening up'' of three already-existing points.

We have found two examples of linear configurations of type $T$ (above) 
whose first infinitesimal neighborhoods have the following two Betti 
diagrams (verified experimentally on {\tt macaulay}). ÊWe are not sure if 
there are any other Betti diagrams possible.

\begin{verbatim}

total: Ê Ê Ê1 Ê Ê 8 Ê Ê 7 Ê Ê Ê Ê total: Ê Ê Ê1 Ê Ê 6 Ê Ê 5
-------------------------- Ê Ê Ê Ê--------------------------
ÊÊ Ê 0: Ê Ê Ê1 Ê Ê - Ê Ê - Ê Ê Ê Ê Ê Ê 0: Ê Ê Ê1 Ê Ê - Ê Ê -
ÊÊ Ê 1: Ê Ê Ê- Ê Ê - Ê Ê - Ê Ê Ê Ê Ê Ê 1: Ê Ê Ê- Ê Ê - Ê Ê -
ÊÊ Ê 2: Ê Ê Ê- Ê Ê - Ê Ê - Ê Ê Ê Ê Ê Ê 2: Ê Ê Ê- Ê Ê - Ê Ê -
ÊÊ Ê 3: Ê Ê Ê- Ê Ê - Ê Ê - Ê Ê Ê Ê Ê Ê 3: Ê Ê Ê- Ê Ê - Ê Ê -
ÊÊ Ê 4: Ê Ê Ê- Ê Ê - Ê Ê - Ê Ê Ê Ê Ê Ê 4: Ê Ê Ê- Ê Ê - Ê Ê -
ÊÊ Ê 5: Ê Ê Ê- Ê Ê - Ê Ê - Ê Ê Ê Ê Ê Ê 5: Ê Ê Ê- Ê Ê - Ê Ê -
ÊÊ Ê 6: Ê Ê Ê- Ê Ê - Ê Ê - Ê Ê Ê Ê Ê Ê 6: Ê Ê Ê- Ê Ê - Ê Ê -
ÊÊ Ê 7: Ê Ê Ê- Ê Ê 4 Ê Ê 2 Ê Ê Ê Ê Ê Ê 7: Ê Ê Ê- Ê Ê 4 Ê Ê -
ÊÊ Ê 8: Ê Ê Ê- Ê Ê 3 Ê Ê 4 Ê Ê Ê Ê Ê Ê 8: Ê Ê Ê- Ê Ê 1 Ê Ê 4
ÊÊ Ê 9: Ê Ê Ê- Ê Ê 1 Ê Ê 1 Ê Ê Ê Ê Ê Ê 9: Ê Ê Ê- Ê Ê 1 Ê Ê 1

\end{verbatim}

\qed
\end{example}

It is possible to isolate an important family of type vectors for which all 
linear configurations of those types have their first infinitesimal 
neighborhoods sharing both the same Hilbert function and same Betti diagram.

\begin{corollary} \label{fat pt unique Hf}
Let $T = (n_1,\dots,n_r)$ be a 2-type vector and let $T' = (m_1,\dots, 
m_{2r})$ be the associated pseudo type vector. ÊIf Ê$n_i \neq 2n_j$ for all 
$i,j$,
then the pseudo type vector $T'$ is actually a 2-type vector. Ê(This holds,
for example, if all the $n_i$ are odd.) Ê In this case Êthe Hilbert 
function and graded Betti numbers of any set of double points supported on 
a linear configuration of type $T$ are uniquely determined, and is that of 
a linear configuration of type $T'$.
\end{corollary}

\begin{proof}
Immediate.
\end{proof}

\section{Beyond linear configurations}

As indicated in the introduction, this paper is intended as a first step in the
study of the following problem: given the Hilbert function $\underline{h}$ 
for a
reduced, zero-dimensional subscheme of $\mathbb P^2$, what are the possible
Hilbert functions of double point schemes whose support has Hilbert function
$\underline{h}$? ÊIn particular, is there a minimum and maximum such function,
$\hu^{\min}, \hu^{\max}$ respectively? ÊIn this section we address these
questions, proving the existence of $\hu^{\max}$ in general and the 
existence of
$\hu^{\min}$ at least in a special case. ÊThe examples in this section also 
help
to clarify the role of linear configurations, and their limitations, toward an
answer to these questions in general.

\begin{example}\label{supp diff hf}
It is not hard to find examples of two sets, $\X$ and $\X'$, of points in
$\mathbb P^2$ with the Êsame Hilbert function, with the property that the
multiplicity two schemes supported on those sets have different Hilbert
functions. Ê A consequence of Theorem \ref{fat = bdl} is that $\X$ and 
$\X'$ can
even have the same graded Betti numbers (e.g.\ both can be linear 
configurations
of the same type), and yet they can have resulting double point schemes with
different Hilbert functions.

A different question is whether there exist unions of double points with the
same Hilbert function, but whose supports have different Hilbert 
functions. ÊThe
answer is ``yes," and we can use Theorem~\ref{any 2-type} to help produce such
an example.

Consider the 2-type vector $(1,2,3,4)$. ÊThis corresponds to a Hilbert function
whose first difference is Ê$\hu_1 = (1,2,3,4)$. ÊThe construction of Theorem
\ref{any 2-type} gives a set, $\Z_1$, of double points whose support, 
$\X_1$, has
Hilbert function with first difference $\hu_1$ and sits on the standard grid,
and such that the Hilbert function of $\Z_1$ has first difference
$(1,2,3,4,5,6,6,3)$.

Now consider a set $\X_2$ of 10 general points on a smooth cubic curve Êin the
plane, Êand let $\Z_2$ be the double points Êsupported on $\X_2$. ÊOne can 
check
with a computer algebra program (e.g.\ {\tt macaulay} \cite{macaulay}) that
$\Z_2$ has the same Hilbert function as described above, and yet $\X_2$ has
Hilbert function with first difference $\hu_2 = (1,2,3,3,1)$. \qed
\end{example}

\begin{example}
It should be noted that this process of studying the Hilbert function of the
first infinitesimal neighborhood of a linear configuration does {\em not} give
all possible Hilbert functions for double points in $\mathbb P^2$. ÊIndeed, a
set of seven generally chosen fat points has Hilbert function whose first
difference is
$(1,2,3,4,5,6)$, while any linear configuration of seven points has at 
least one
subset of four points on a line, so the regularity must be at least 8 for the
corresponding double points.\qed
\end{example}

We now turn to the question of the existence of $\hu^{\max}$ and $\hu^{\min}$.
We are grateful to Mike Roth for useful discussions about the following theorem
and its proof.

\begin{theorem} \label{existence of hmax}
Let $\hu$ be the Hilbert function of some reduced zero-dimensional subscheme of
$\mathbb P^2$. ÊThen there is a Hilbert function $\hu^{\max}$ such that if $h'$
is the Hilbert function of a double point scheme whose support has Hilbert
function $\hu$ then $h' \leq \hu^{\max}$.
\end{theorem}

\begin{proof}
Let $Hilb^s(\P^2) = \X^{(s)}$ be the Hilbert scheme which parametrizes all 
the closed subschemes of $\P^2$ having length $s$. ÊIt is well known (e.g. 
by using the Hilbert-Burch Theorem) that those closed subschemes of $\P^2$ 
which share the same Hilbert function, $\underline{g}$ (say), form an 
irreducible subset of $\X^{(s)}$ (which we'll denote by 
$\X^{(s)}_{\underline{g}}$). ÊIt is also well known that 
$\X^{(s)}_{\underline{g}}$ is locally closed. ÊThus, for any positive 
integer $s$ we obtain a (finite) locally closed irreducible partition of 
$Hilb^s(\P^2)$.

The partitiion we described above gives, in the same way, a partition of 
$Sym^s(\P^2) = \Y^{(s)}$ (the scheme parametrizing families of $s$ distinct 
points in $\P^2$).

There is also a map from $\Y^{(t)} = Sym^t(\P^2)$ into $\X^{(3t)} = 
Hilb^{3t}(\P^2)$ -- we associate to a set of $t$ distinct points in $\P^2$ 
its first infinitesimal neighborhood. ÊWe'll denote the image of $\Y^{(t)}$ 
in $Hilb^{3t}(\P^2)$ by ${\mathcal D}^{(t)}$ and the image of 
$\Y^{(t)}_{\underline{h}}$ by ${\mathcal D}^{(t)}_{\underline{h}}$.

If we restrict the stratification of $Hilb^{3t}(\P^2)$ to ${\mathcal 
D}^{(t)}_{\underline{h}}$ then exactly one component of this stratification 
will be dense in ${\mathcal D}^{(t)}_{\underline{h}}$, and this stratum will be
$$
{\mathcal D}^{(t)}_{\underline{h}} \cap {\X}^{(3t)}_{\underline{g}}
$$
for some $\underline{g}$. ÊThat $\underline{g} = \hu^{\max}$.
\end{proof}

\begin{remark}
Theorem \ref{existence of hmax} is an existence result, valid for any 
Hilbert function $\hu$. ÊUnfortunately, we do not know (in general) an 
explicit formula (or even an algorithm) for computing it. ÊHowever, in 
certain special cases we can given an algorithm that easily leads to 
$\hu^{\max}$.

First, suppose that $\hu$ is the Hilbert function of a complete 
intersection of type $(a,b)$. ÊThen, in the irreducible family of sets of 
points with Hilbert function $\hu$, an open subset corresponds to the 
complete intersections of Êtype $(a,b)$. ÊBut, it is well known that if $I 
= (F,G)$ is the ideal of such a complete intersection, $\X$, then $I^2 = 
(F^2, G^2):I$. ÊSince, for a complete intersection $\X$, $I^2$ is the 
defining ideal of the first infinitesimal neighborhood of $\X$, easy 
Liaison techniques give the Hilbert function of $I^2$, which is $\hu^{\max}$.

Second, suppose that $\hu$ corresponds to the 2-type vector $(n_1, n_2, 
\ldots , n_r)$ with $n_i \geq n_{i-1} + 3$ for all $i \geq 2$. ÊThen any 
reduced set of points whose Hilbert function has this type vector must be a 
$k$-configuration (using the decompostion techniques of Davis 
\cite{davis}). ÊThe general such $k$-configuration is a linear 
configuration and the Hilbert function of its first infinitesimal 
neighborhood is uniquely determined by Theorem \ref{fat = bdl}. ÊTherefore, 
this is $\hu^{\max}$.
\end{remark}

ÊÊWe have been unable to prove that $\hu^{\min}$ exists, in 
general. ÊHowever, we will prove its existence in an important special 
case, and give a conjecture for the general case. ÊIn what follows we 
continue our abuse of notation and refer to a curve and its defining form 
interchangably.

\begin{lemma} \label{max dble}
Let $F$ be a reduced curve of degree $d$.

\begin{itemize}

\item[(a)] If $F$ is a union of $d$ lines, each of which meets the
remaining lines in $d-1$ distinct points, then the number of singular
points of $F$ is $\binom{d}{2}$, all double points.

\item[(b)] If $F$ is not a union of $d$ lines, each of which meets the
remaining lines in $d-1$ distinct points, then the number of singular
points of $F$ is $< \binom{d}{2}$.
\end{itemize}
\end{lemma}

\begin{proof}
(a) is clear. ÊFor (b), suppose first that $F$ is irreducible. ÊThen the
number of double points is $\leq \frac{(d-1)(d-2)}{2} < \binom{d}{2}$.
Now suppose that $F$ is not irreducible, $F = F_1 \cdot F_2$. ÊIf $F_1$
and $F_2$ are both unions of lines but at least three lines pass through
one point then clearly the number of singular points is $<
\binom{d}{2}$. ÊFinally, suppose that $F = F_1\cdot F_2$ where at least
one, say $F_1$, is irreducible of degree $\geq 2$. ÊSay $\deg F_i = d_i$,
with $d_1 + d_2 = d$. ÊBy induction, then, the number of singular points
of $F_1$ is $< \binom{d_1}{2}$ while the number of singular points of
$F_2$ is $\leq \binom{d_2}{2}$, with equality if and only if $F_2$ is a
suitable union of lines. ÊThe singular points of $F$ then come either
as singular points of $F_1$ or $F_2$, or as points of intersection of
$F_1$ and $F_2$. ÊThen the number of singular points of
$F$ is
\[
\# Sing(F) < \binom{d_1}{2} + \binom{d_2}{2} + d_1 d_2 =
\binom{d_1+d_2}{2} = \binom{d}{2}
\]
as desired.
\end{proof}

\begin{notation} \label{def of Ctr}
Let $\lambda_1,\dots,\lambda_t$ be a set of $t$ distinct lines in
$\mathbb P^2$ such that each $\lambda_i$ meets the remaining $t-1$ lines
in $t-1$ distinct points. ÊWe denote by $C_t$ the configuration
consisting of the $\binom{t}{2}$ pairwise intersections of these lines.
Let $0 \leq r \leq t$. ÊWe denote by $C_{t,r}$ a subconfiguration of
$C_{t+1}$ obtained by removing any $(t-r)$ points of $C_{t+1}$ that lie
on $\lambda_{t+1}$. ÊNote that $C_t \subseteq C_{t,r} \subseteq C_{t+1}$.
The first equality holds if $r=0$ and the second holds if $r=t$. \qed
\end{notation}

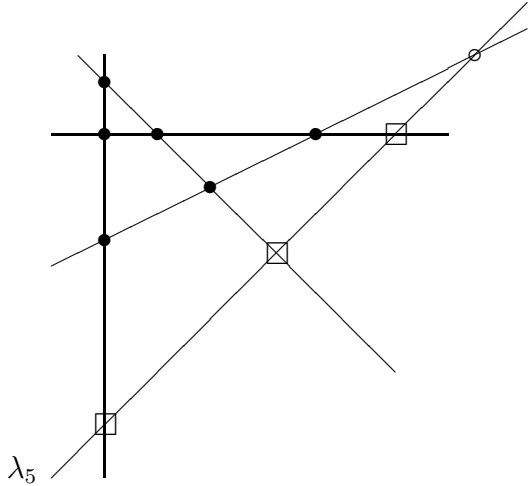
\begin{example}

\vbox{$$
\begin{picture}(130,140)(-20,-100)
\put(-50,-20){\line(1,0){150}}
\put(-40,10){\line(1,-1){120}}
\put(-30,10){\line(0,-1){160}}
\put(-50,-70){\line(2,1){180}}
\put(-50,-150){\line(1,1){180}}
\put(-33,-3){$\bullet$}
\put(-33,-23){$\bullet$}
\put(-33,-63){$\bullet$}
\put(-46,-134){$\qed$}
\put(-13,-23){$\bullet$}
\put(47,-23){$\bullet$}
\put(64,-24){$\qed$}
\put(7,-43){$\bullet$}
\put(19,-69){$\qed$}
\put(107,7){$\circ$}
\put(-67,-150){$\lambda_5$}
\end{picture}$$}

\vskip .9in

The bullets represent $C_4$. ÊThe bullets together with the squares
represent $C_{4,3}$. ÊThe bullets, squares and circle together represent
$C_5$. ÊNote that $C_4 \subset C_{4,3} \subset C_5$. \qed
\end{example}

\begin{lemma} \label{details Ctr}
\hbox{\hskip 1in}

\begin{itemize}
\item[(a)] $\deg C_t = \binom{t}{2}$

\item[(b)] $\deg C_{t,r} = \binom{t}{2} + r$

\item[(c)] The first difference of the Hilbert function of $C_t$ is
\[
\begin{array}{cccccccccccccccccccccccc}
1 & 2 & 3 & \dots & (t-1)
\end{array}
\]
ÊÊand the first difference of the Hilbert function of $C_{t,r}$ is
\[
\begin{array}{ccccccccccccccccccccc}
1 & 2 & 3 & \dots & (t-1) & r
\end{array}
\]
In particular, $C_{t,r}$ has so-called {\em generic Hilbert function}.
\end{itemize}
\end{lemma}

\begin{proof}
(a) and (b) are clear. ÊFor the first part of (c), suppose that $C_t$
lies on a curve $F$ of degree $t-2$. ÊEach line $\lambda_i$ contains
$t-1$ collinear points of $C_t$, so by Bezout's theorem $\lambda_i$ must be
a component of $F$. ÊBut there are $t$ such lines. ÊContradiction. ÊThe
second part of (c) comes from the first part together with the inclusions
$C_t \subset C_{t,r} \subset C_{t+1}$, and the fact that consequently the
first difference of the Hilbert function of $C_{t,r}$ must be between those of
$C_t$ and
$C_{t+1}$.
\end{proof}

\begin{notation}
We denote by $\Z_t$ the first infinitesimal neighborhood of $C_t$. Ê We
denote by $\Z_{t,r}$ the first infinitesimal neighborhood of $C_{t,r}$.
Note that $\Z_t = \Z_{t,0}$. \qed
\end{notation}

\begin{theorem} \label{unique Zt}
\begin{itemize}
\item[(a)] The first difference of the Hilbert function of $\Z_t$ is

\begin{center}
\begin{tabular}{r|cccccccccccccccccc}
\hbox{\rm degree} & $0$ & $1$ & $2$ & $3$ & Ê$\dots$ & $(t-1)$ & $t$ &
$(t+1)$ &
$
\dots
$ &
$2t-3$ & $2t-2$ \\ \hline
$\Delta h_{\Z_t}$ & $1$ & $2$ & $3$ & $4$ &$\dots$ & $t$ & $t$ & $t$ &
$\dots$ &
$t$ &
$0$
\end{tabular}

\end{center}

ÊÊ\medskip

Note that there are $t-1$ occurrences of $t$ at the end of this
function.

\item[(b)] Among double point schemes whose support has Hilbert function with
first difference $\hu = (1, 2, 3 , \dots , t-1)$, $\Z_t$ has minimal Hilbert
function.

\item[(c)] Up to a different choice of $\lambda_1,\lambda_2, \dots ,
\lambda_t$, $\Z_t$ is the unique double point scheme with this Hilbert
function, among double point schemes whose support has Hilbert function with
first difference $\hu$. ÊIn fact, the value of this Hilbert function in degree
$t$ already uniquely determines $\Z_t$.
\end{itemize}
\end{theorem}

\begin{proof}
For (a), first note that the union of the lines
$\lambda_1,\dots,\lambda_t$ is a component of any curve of degree $\leq
2(t-1)-1 = 2t-3$ containing $\Z_t$, by Bezout's theorem. ÊOn the other
hand, this union is double at each of the $\binom{t}{2}$ points of
$C_t$. ÊHence the ideal has exactly one generator in degree $t$, and the
next generator does not come before degree $2t-2$. ÊSo the first difference of
the HIlbert function of
$\Z_t$ must be as claimed at least up to degree $2t-3$. ÊBut
\[
1+2+\dots + (t-1)+t + t + \dots + t = \binom{t}{2} + (t-1)t = 3
\binom{t}{2} = \deg \Z_t,
\]
so this must be the full Hilbert function.

We now prove (b) and (c) at the same time. ÊLet $\X$ be a reduced set of
$\binom{t}{2}$ points with generic Hilbert function (i.e.\ the one with
Hilbert function with first difference as given in Lemma \ref{details Ctr} (c)
for $C_t$) and let $\Z$ be its first infinitesimal neighborhood. ÊSuppose that
$\Z$ has Hilbert function that is strictly smaller than that of $\Z_t$ in some
degree. ÊWe consider the Êfirst difference of the Hilbert function of $\Z$, 
first
in degree
$t-1$. ÊSuppose that
$h_\Z(t-1) < h_{\Z_t}(t-1)$, i.e.\ suppose that $\Z$ lies on some curve
$F$ of degree $t-1$. ÊThen $F$ is at least double at all the points of
$\X$. ÊSince $\X$ has the generic Hilbert function of Lemma \ref{details
Ctr} (c), the initial degree of $I_\X$ is
$t-1$, and in particular it lies on no curve of degree $t-2$. ÊBy
assumption there is a form $F$ of degree $t-1$ containing $\X$ that is in
fact (at least) double at all the points of $\X$. ÊWe first claim that
$F$ is reduced. ÊIf it were not, then the radical is a form of degree
$<t-1$ containing $\X$, contradicting the fact that $t-1$ is the initial
degree of $I_\X$. ÊBut now Lemma \ref{max dble} says that $F$ has at most
$\binom{t-1}{2}$ singular points. ÊThis contradiction shows that
$h_\Z(t-1) = h_{\Z_t}(t-1)$.

We now turn to degree $t$. ÊSuppose that the initial degree of $I_\Z$ is
$t$, so $h_\Z(t) \leq h_{\Z_t}(t)$. ÊThen there is at least one form, $F$,
of degree $t$ that is singular at all the points of $\X$.

\bigskip

\noindent {\bf Claim:} $F$ is reduced.

\medskip

To prove this claim, we suppose otherwise. ÊThen $F$ has a factor, $F_1$,
that is not reduced. ÊIf $\deg F_1 \geq 2$ then the radical of $F$ is a
form of degree $\leq t-2$ that contains $\X$, again contradicting the
fact that the initial degree of $I_\X$ is $t-1$. ÊSo now suppose that
there is a linear form, $L$, such that $F = L^2 F_2$, with $F_2$
reduced. ÊFrom the first difference of the Hilbert function of $\X$, we see 
that
$\X$ contains at most
$t-1$ collinear points. ÊHence $F_2$ is a reduced form of degree $t-2$
double at $\binom{t}{2} - (t-1) = \binom{t-1}{2}$ points or more. ÊThis
violates Lemma \ref{max dble} and proves our Claim.

\bigskip

So now we have $h_\Z(t) \leq Êh_{\Z_t}(t)$ and $I_\Z$ contains a reduced
form, $F$, of degree $t$ that is double at $\binom{t}{2}$ points. ÊBy
Lemma \ref{max dble}, then, $F$ is a union of lines and $\X$ is the
pairwise intersection of these lines. ÊSo $\Z = \Z_t$ (up to the choice
of $\lambda_i$).

We may thus assume without loss of generality that the initial degree of
$I_\Z$ is $\geq t+1$, so the first difference of the Hilbert function of 
$\Z$ is

\medskip

\begin{center}
\begin{tabular}{r|cccccccccccccccccc}
\hbox{\rm degree} & $0$ & $1$ & $2$ & $3$ & Ê$\dots$ & $(t-2)$ & $(t-1)$ &
$t$ & $(t+1)$ &
$
\dots
$ Ê\\ \hline
$\Delta h_{\Z}$ & $1$ & $2$ & $3$ & $4$ & $\dots$ & $t-1$ & $t $ & $t+1$ &
? &
$\dots$
\end{tabular}

\end{center}

\medskip

Recall that the first difference of the Hilbert function of $\Z_t$ is

\begin{center}
\begin{tabular}{r|cccccccccccccccccc}
\hbox{\rm degree} & $0$ & $1$ & $2$ & $3$ & Ê$\dots$ & $(t-2)$ & $(t-1)$ &
$t$ & $(t+1)$ &
$
\dots
$ &
$2t-3$ & $2t-2$ \\ \hline
$\Delta h_{\Z_t}$ & $1$ & $2$ & $3$ & $4$ &$\dots$ &$(t-1)$ & $t$ & $t$ &
$t$ &
$\dots$ &
$t$ &
$0$
\end{tabular}

\end{center}

\medskip

\noindent In particular, we have $h_\Z(t) > h_{\Z_t}(t)$.
We have to show that it cannot happen that later on, the Hilbert function
of $\Z$ drops below that of ${\Z_t}$. ÊBy Lemma \ref{bd reg}, the
regularity of $I_\Z$ is $\leq 2 \cdot \reg(I_\X) = 2(t-1) = 2t-2$, so the
first difference of the Hilbert function of $\Z$ ends in degree $\leq 2t-3$ as
well.

Suppose that there is a value, $d$, for which $h_\Z(d) < h_{\Z_t}(d)$.
Clearly $t < d < 2t-3$, since $h_\Z(2t-3) = h_{\Z_t}(2t-3) = 3
\binom{t}{2}$.
The Hilbert function in any degree is just the sum of the entries of the
first difference, up to and including that degree. ÊBut since $h_\Z(t) >
h_{\Z_t}(t)$, this means that the first difference of the Hilbert function of
$\Z$ in some degree
$\leq d$ has a value $k<t$. ÊBut the first difference of the Hilbert 
function of
a zero-dimensional subscheme of $\mathbb P^2$ is non-increasing in degrees 
$\geq \alpha$ (see Definition \ref{hilb} $i)$\ ). ÊHence
\[
\begin{array}{rcl}
\deg \Z_t & = & h_{\Z_t}(d) + t(2t-3-d) \\
& > & h_\Z(d) + t(2t-3-d) \\
& > & h_\Z(d) + k(2t-3-d) \\
& \geq & \deg \Z.
\end{array}
\]
This contradiction shows that $\Z_t$ does in fact have minimal Hilbert
function as claimed.
\end{proof}

\begin{remark} \label{good example}
Theorem \ref{unique Zt} illustrates the necessity of restricting our 
hypothesis in Theorem \ref{fat = bdl} to linear configurations for the 
support rather than $k$-configurations. Ê First note that $C_t$ is a 
$k$-configuration but not a linear configuration. ÊIndeed, every newly 
added line misses all previous points of the configuration, but the points 
on the new line do lie on previously
existing lines.

We now consider an example. ÊLet $t=4$. ÊThen the Êfirst difference of the
Hilbert function of
$\Z_4$ is

\medskip

\begin{center}
\begin{tabular}{r|cccccccccccccccccc}
\hbox{\rm degree} & $0$ & $1$ & $2$ & $3$ & 4 & 5 & 6 \\ \hline
$\Delta h_{\Z_4}$ & $1$ & $2$ & $3$ & $4$ & 4 & 4 & 0
\end{tabular}

\end{center}

\medskip

\noindent On the other hand, the configuration $C_4$ has Hilbert function with
first difference
$(1,2,3)$, so this is also the type vector (in this case). ÊThe associated
pseudo type vector is
$(1,2,2,3,4,6)$. ÊBy Theorem \ref{fat = bdl}, however, any linear configuration
with type vector $(1,2,3)$ has first infinitesimal neighborhood whose Hilbert
function has first difference
$(1,2,3,4,5,3)$.

This example also serves as a counterexample to a natural guess, namely
that the standard configuration (or in general the spread out
configuration) should yield the first infinitesmial neighborhood of
minimal Hilbert function among all supports with fixed Hilbert function.
Indeed, the problem is that these configurations $C_t$ have even more
collinearities than the spread out configurations. \qed

\end{remark}

We now consider generic Hilbert functions $\underline{h}$ that do not
correspond to precisely $\binom{t}{2}$ points. One would like to find the
minimal Hilbert function, $\underline{h}^{\min}$, for the first
infinitesmial neighborhoods of point sets with Hilbert function
$\underline{h}$.
For example, we now compute the first difference of the Hilbert functions of
some low-degree examples.

\medskip

\begin{center}

\begin{tabular}{l|cccccccccccccccccc}
\hbox{\rm degree} & $0$ & $1$ & $2$ & $3$ & 4 & 5 & 6 & 7 & 8 & 9 & 10 \\
\hline
$\Delta h_{\Z_4}$ & $1$ & $2$ & $3$ & $4$ & 4 & 4 & Ê\\
$\Delta h_{\Z_{4,1}}$ & 1 & 2 & 3 & 4 & 5 & 4 & 1 & 1 \\
$\Delta h_{\Z_{4,2}}$ & 1 & 2 & 3 & 4 & 5 & 5 & 2 & 2 \\
$\Delta h_{\Z_{4,3}}$ & 1 & 2 & 3 & 4 & 5 & 5 & 4 & 3 \\
$\Delta h_{\Z_5}$ & 1 & 2 & 3 & 4 & 5 & 5 & 5 & 5 \\
$\Delta h_{\Z_{5,1}}$ & 1 & 2 & 3 & 4 & 5 & 6 & 5 & 5 & 1 & 1 \\
$\Delta h_{\Z_{5,2}}$ & 1 & 2 & 3 & 4 & 5 & 6 & 6 & 5 & 2 & 2 \\
$\Delta h_{\Z_{5,3}}$ & 1 & 2 & 3 & 4 & 5 & 6 & 6 & 6 & 3 & 3 \\
$\Delta h_{\Z_{5,4}}$ & 1 & 2 & 3 & 4 & 5 & 6 & 6 & 6 & 5 & 4 \\
$\Delta h_{\Z_6}$ & 1 & 2 & 3 & 4 & 5 & 6 & 6 & 6 & 6 & 6
\end{tabular}

\end{center}

\medskip

For instance, why should $h_{\Z_{4,2}}$ be minimal?

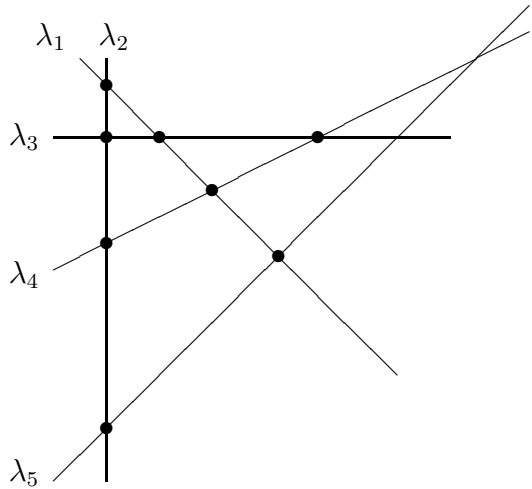
\vbox{$$
\begin{picture}(130,140)(-20,-100)
\put(-50,-20){\line(1,0){150}}
\put(-40,10){\line(1,-1){120}}
\put(-30,10){\line(0,-1){160}}
\put(-50,-70){\line(2,1){180}}
\put(-50,-150){\line(1,1){180}}
\put(-33,-3){$\bullet$}
\put(-33,-23){$\bullet$}
\put(-33,-63){$\bullet$}
\put(-33,-133){$\bullet$}
\put(-13,-23){$\bullet$}
\put(47,-23){$\bullet$}
\put(7,-43){$\bullet$}
\put(32,-68){$\bullet$}
\put(-67,-150){$\lambda_5$}
\put(-67, -75){$\lambda_4$}
\put(-67, -23){$\lambda_3$}
\put(-57, 15){$\lambda_1$}
\put(-33,15){$\lambda_2$}
\end{picture}$$}

\vskip .9in

\noindent We have not been able to find an argument even in this 
case. ÊHowever, we have the following:

\begin{conjecture}
Among double schemes whose support has a fixed generic Hilbert function
$(1,2,\dots, t-1, r)$ (see Lemma \ref{details Ctr}), there is a minimal Hilbert
function, and it occurs when the support is $C_{t,r}$.
\end{conjecture}

Note that when $0 < r < t$, we do not conjecture that the minimal Hilbert
function can {\em only} occur when the support is $C_{t,r}$, as was the case
for $C_t$. ÊFor instance, the Hilbert function for the first infinitesimal
neighborhood of $C_{4,2}$ (illustrated above) can also arise from the first
infinitesimal neighborhood of the following configuration (where the oval is a
conic):

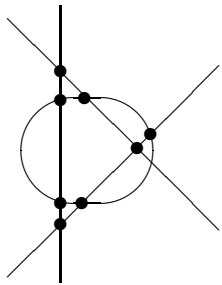
\vbox{$$
\begin{picture}(130,140)(-20,-100)
\put(-50,10){\line(1,-1){80}}
\put(-30,15){\line(0,-1){105}}
\put(-50,-88){\line(1,1){80}}
\put(-20,-40){\oval(50,40)}
\put(-33,-13){$\bullet$}
\put(-33,-24){$\bullet$}
\put(-33,-63){$\bullet$}
\put(-33,-71){$\bullet$}
\put(-24,-23){$\bullet$}
\put(-25,-63){$\bullet$}
\put(1,-37){$\bullet$}
\put(-4,-42){$\bullet$}
\end{picture}$$}

\medskip\medskip

We generalize the above to any Hilbert function corresponding to type 
$(n_1, \ldots , n_r)$. ÊWe form a configuration $C_{\hu}$ as 
follows. ÊChose a set of $r+1$ lines (we'll call them $\lambda _1, \ldots , 
\lambda _{r+1}$) and let $C_{r+1}$ be as above, i.e. the union of all the 
pairwise intersection points of the $\lambda _i$'s. ÊSo, each of the 
$\lambda_i$ contains $r$ points of $C_{r+1}$. ÊNotice, however, that we can 
view $C_{r+1}$ as a $k$-configuration in the following way: first choose 
all $r$ points on $\lambda _r$; then, on $\lambda_{r-1}$, choose the 
remaining $r-1$ points (since one was already chosen on $\lambda _r$); on 
$\lambda _{r-2}$ choose the remaining $r-2$ points; \ldots ; on $\lambda 
_1$ choose the only point remaining. Ê(Note that $\lambda_{r+1}$ has become 
irrelevant in this point of view.) ÊNow we add $n_r - r$ arbitrary points 
on $\lambda _r$, $n_{r-1}-(r-1)$ points on $\lambda _{r-1}$ , etc. thus 
forming a $k$-configuration, $C_{\hu}$ of type $(n_1, \ldots , n_r)$.

\medskip\noindent\begin{conjecture}
The Hilbert function of $C_{\hu}$ is $\hu^{\min}$.
\end{conjecture}



\begin{thebibliography}{999}

\bibitem{AH3} J.\ Alexander and A.\ Hirschowitz, {\em Polynomial interpolation
in several variables}, J.\ Alg.\ Geom.\ {\bf 4} (1995), 201--222.

\bibitem{macaulay} D.\ Bayer and M.\ Stillman, Macaulay: A system for
computation in algebraic geometry and commutative algebra. Source and object
code available for Unix and Macintosh computers. ÊContact the authors, or
download from ftp://math.harvard.edu via anonymous ftp.

\bibitem{bazz} L.\ Bazzotti, {\em ??? Conductors for schemes in $\P^2$}

\bibitem{bigatti} A.\ Bigatti, {\em Upper bounds for the Betti numbers of a
given Hilbert function}, Comm.\ Algebra {\bf 21} (1993), no. 7, 2317--2334.

\bibitem{BGM} A.\ Bigatti, A.V.\ Geramita and J.\ Migliore, {\em Geometric
Consequences of Extremal Behavior in a Theorem of Macaulay}, Trans.\ ÊAmer.\
Math.\ Soc.\ {\bf 346} (1994), 203--235.

\bibitem{BM4} G.\ Bolondi and J.\ Migliore, {\em The Structure of an Even
Liaison ÊClass}, Trans.\ Amer.\ Math.\ Soc.\ {\bf 316} (1989), 1--37.

\bibitem{bocci} C.\ Bocci, {\em Special linear systems and special effect 
varieties}, Ph.D. Thesis, University of Torino (2003).

\bibitem{bruns-herzog} W.\ Bruns and J.\ Herzog, ``Cohen-Macaulay Rings,''
Cambridge studies in advanced mathematics, Cambridge University Press, 1993.

\bibitem{campanella} G.\ Campanella, {\em Standard bases of perfect 
homogeneous polynomial ideals of height 2}, J.\ Alg. \ {\bf 101}, (1986), 
47-60.

\bibitem{catalisano} M.V.\ Catalisano, {\em ``Fat" points on a conic}, 
Comm.\ in Alg. \ {\bf 19} (1991) 2153-2168.

\bibitem{cgg} M.V.\ Catalisano, A.V.\ Geramita, A.\ Gimigliano, {\em Higher 
Secant Varieties of Segre-Veronese Varieties}, {\it to appear} Proceedings 
of the Conference in Honour of G. Veronese, Siena (2004).

\bibitem{cgg2} M.V.\ Catalisano, A.V.\ Geramita, A.\ Gimigliano {\em Higher 
Secant Varieties of the Segre Varieties $\P^1 \times \cdots \times \P^1$ } 
To appear, Jour. of Pure and Appl. Algebra.

\bibitem{cat-el-gim} M.V.\ Catalisano, P.\ Ellia, A.\ Gimigliano {\em Fat 
points on rational normal curves} Jo.\ Alg.\ {\bf 216} (1999) 600-619.

\bibitem{chandler} K. \ Chandler, {\em A brief proof of a maximal rank 
theorem for generic double points in projective space}, Trans.\ Amer.\ 
Math.\ Soc. \ {\bf 353} (2001), no.\ 5, 1907-1920.

\bibitem{ciliberto} C.\ Ciliberto, {\em Geometric aspects of polynomial 
interpolation in more variables and of Waring's problem} \ Proceedings of 
the Third European Congress of Mathematics, I (Barcelona, 2000), 289-316, 
Prog.\ Math.\ {\bf 201} Birkhauser, Basel, 2001.

\bibitem{ciliberto-miranda} C.\ Ciliberto and R.\ Miranda, {\em 
Degenerations of
planar linear systems}, J.\ Reine Angew.\ Math.\ {\bf 501} (1998), 191--220.

\bibitem{ciliberto-miranda2} C.\ Ciliberto and R.\ Miranda, {\em Linear 
systems of
plane curves with base points of equal multiplicity}, Trans.\ Amer.\ Math.\ 
Soc.\ {\bf
352} (2000), no.\ 9, 4037--4050.

\bibitem{davis} E.D.\ Davis, {\em Complete Intersections of Codimension 2 in
$\mathbb P^r$: The Bezout-Jacobi-Segre Theorem Revisited}, Rend.\ Sem.\ Mat.\
Univers.\ Politecn.\ Torino, Vol.\ 43, 4 (1985), 333--353.

\bibitem{dent} A.\ Dent, {\em Variations on Methods of Lorentz and Lorentz 
for Dimensions Two and Three,} Ph.D. Thesis, Colorado State University (2003)

\bibitem{diesel} S.\ J.\ Diesel, {\em Irreducibility and dimension theorems 
for families of height 3 Gorenstein algebras}, Pac. J. Math. {\bf 172} 
(1996) 365-397.

\bibitem{fat} G.\ Fatabbi, {\em On the resolution of ideals of fat points}, 
J. Alg.,{\bf 242} (2001), 92-108.

\bibitem{fat-lor} G.\ Fatabbi, A.\ Lorenzini {\em On the graded resolution 
of ideals of a few general fat points of $\P^n$,} To appear, J. Pure Appl. 
Algebra.

\bibitem{FHH} S.\ Fitchett, B.\ Harbourne, S.\ Holay {\em Resolutions of 
Fat Point Ideals involving 8 General Points of $\P^2$} J.\ of\ Alg.\ {\bf 
24} (2001) 684-705.

\bibitem{fran} C.A.\ Francisco, {\em Resolutions of small sets of fat 
points} (preprint, 2004).

\bibitem{ghs} A.V.\ Geramita, T.\ Harima and Y.S.\ Shin, {\em An Alternative
to the Hilbert Function for the Ideal of a Finite Set of Points in
$\mathbb P^n$}, Illinois J.\ Math.\ {\bf 45} (2001), no. 1, 1--23..

\bibitem{GHSAdv} A.V.\ Geramita, T.\ Harima and Y.S.\ Shin, {\em Extremal Point
Sets and Gorenstein Ideals}, ÊAdv.\ Math.\ {\bf 152} (2000), 78--119.

\bibitem{GMR} A.V.\ Geramita, P.\ Maroscia and L.\ Roberts, {\em The Hilbert
Function of a Reduced $k$-Algebra}, J.\ London Math.\ Soc.\ {\bf 28} (1983),
443--452.

\bibitem{GM4} A.V.\ Geramita and J.\ Migliore, {\em A Generalized Liaison
Addition}, J.\ Alg.\ 163 (1994), 139--164.

\bibitem{GM} A.V.\ Geramita and J.\ Migliore, {\em Reduced Gorenstein 
Codimension 3 Subschemes of Projective Space}, Proc. A.M.S.,{\bf 125} 
(1997), 943-950.

\bibitem{GS} A.V.\ Geramita and Y.S.\ Shin, {\em $k$-configurations in
$\mathbb P^{3}$ all have Extremal Resolutions}, ÊJ.\ ÊAlgebra, {\bf 213} 
(1999),
351--368.

\bibitem{harbourne} B.\ Harbourne, {\em The geometry of rational surfaces and
Hilbert functions of points in the plane}, Proc.\ 1984 Vancouver Conf.\ on 
Alg.\
Geom.\, A.M.S., Providence, RI 1986, 95--111.

\bibitem{HHF} B.\ Harbourne, S.\ Holay, S.\ Fitchett, {\em Resolutions of 
Ideals of Quasiuniform Fat Point Subschemes of $\P^2$.} Trans. Amer. Math. 
Soc. {\bf 355} (2002) 593-208.

\bibitem{hulett} H.\ Hulett, {\em Maximum Betti numbers of homogeneous ideals
with a given Hilbert function}, Comm.\ Algebra {\bf 21} (1993), no. 7,
2335--2350.

\bibitem{KMMNP} J.\ Kleppe, J.\ Migliore, R.M.\ Mir\'o-Roig, U.\ Nagel and C.\
Peterson, {\em Gorenstein Liaison, Complete Intersection Liaison Invariants and
Unobstructedness}, Memoirs of the Amer.\ Math.\ Soc.\ Vol.\ 154, 2001; 116 pp.
Softcover, ISBN 0-8218-2738-3.

\bibitem{laface-ugaglia} ÊA.\ Laface and L.\ Ugaglia, {\em On standard
transformations of $\mathbb P^n$ and special linear series}, preprint {\tt
mathAG/0409129}.

\bibitem{laface-ugaglia2} A.\ Laface and L.\ Ugaglia, {\em 
Quasi-homogeneous linear systems on $\P^2$ with base points of multiplicity 
5}, preprint {\tt
mathAG/0205270}.

\bibitem{LR} R.\ Lazarsfeld and P.\ Rao, {\em Linkage of General Curves of 
Large
Degree}, in ``Algebraic Geometry-- Open Problems (Ravello, 1982),'' Lecture
Notes in Mathematics, vol.\ 997, Springer--Verlag (1983), 267--289.


\bibitem{mig-book} ÊJ.\ Migliore, ``Introduction to Liaison Theory and
Deficiency Modules,'' ÊBirkh\"auser, Progress in Mathematics 165, 1998; 224
pp. Hardcover, ISBN 0-8176-4027-4.


\bibitem{MN1} J.\ Migliore and U.\ Nagel, {\em On the Cohen-Macaulay Type 
of the
General Hypersurface Section of a Curve}, ÊMath.\ Zeit.\ {\bf 219 (2)} (1995),
245--273.

\bibitem{MN2} J.\ Migliore and U.\ Nagel {\em Reduced Arithmetically 
Gorenstein Schemes and Simplicial Polytopes with Maximal Betti Numbers}, 
Adv.\ Math.\ {\bf 180} (2003) 1-63.

\bibitem{miranda} R.\ Miranda, {\em Linear Systems of Plane Curves}, 
Notices A.M.S., {\bf 46}, (1999), 192-202.

\bibitem{mumford} D.\ Mumford, Ê``Lectures on curves on an algebraic surface.
With a section by G. M. Bergman.'' Annals of Mathematics Studies, No.\ 59 
Princeton
University Press, Princeton, N.J.\ 1966 xi+200~pp.

\bibitem{sindipaper} L.\ Sabourin, {\em $n$-type vectors and the 
Cayley-Bacharach
property}, Comm.\ Algebra {\bf 30} (2002), no.8, 3891-3915.

\bibitem{pardue} K.\ Pardue, {\em Deformation Classes of Graded Modules and 
Maximal Betti Numbers}, Illinois J. of Math.\ \textbf{40} (1996), 564--585.

\bibitem{valla} G.\ Valla, {\em Betti numbers of some monomial ideals}, 
Proc. Amer. Math. Soc. {\bf 133} (2005) , 57-63.

\bibitem{yang} S.\ Yang, {\em Linear series in $\mathbb P^2$ with base points
of bounded multiplicity}, preprint {\tt mathAG/0406591}.

\end{thebibliography}
\end{document}